\documentclass[11pt]{article}
\usepackage[utf8]{inputenc}
\usepackage{amsmath,amssymb,epsfig,bbm}
\usepackage{stmaryrd}
\usepackage{comment}
\usepackage[dvips]{color}
\usepackage[T1]{fontenc}

\usepackage{todonotes}
\usepackage{enumitem}
\setlist{nolistsep}

\newtheorem{defi}{Definition}
\newtheorem{prop}[defi]{Proposition}
\newtheorem{theo}[defi]{Theorem}
\newtheorem{conj}[defi]{Conjecture}
\newtheorem{lemm}[defi]{Lemma}
\newtheorem{coro}[defi]{Corollary}
\newtheorem{rema}[defi]{Remark}
\newtheorem{exem}[defi]{Example}
\newtheorem{exems}[defi]{Examples}

\newcommand{\bdefi}{\begin{defi}}
\newcommand{\edefi}{\end{defi}}
\newcommand{\bprop}{\begin{prop}}
\newcommand{\eprop}{\end{prop}}
\newcommand{\btheo}{\begin{theo}}
\newcommand{\etheo}{\end{theo}}
\newcommand{\blemm}{\begin{lemm}}
\newcommand{\brema}{\begin{rema}}
\newcommand{\erema}{\end{rema}}
\newcommand{\bexer}{\begin{exem}}
\newcommand{\eexer}{\end{exem}}
\newcommand{\bexems}{\begin{exems}}
\newcommand{\eexems}{\end{exems}}
\newcommand{\bconj}{\begin{conj}}
\newcommand{\econj}{\end{conj}}
\newcommand{\elemm}{\end{lemm}}
\newcommand{\bcoro}{\begin{coro}}
\newcommand{\ecoro}{\end{coro}}
\newcommand{\dem}{\noindent{\bf Proof. }}
\newcommand{\rem}{\noindent{\bf Remark. }}

\usepackage{mathrsfs}
\renewcommand\mathcal{\mathscr}

\newcommand{\N}{{\cal N}}

\renewcommand{\H}{{\cal H}}

\newcommand{\C}{{\cal C}}

\newcommand{\maths}[1]{{\mathbb #1}}  

\newcommand{\RR}{\maths{R}}
\newcommand{\NN}{\maths{N}}
\newcommand{\CC}{\maths{C}}

\newcommand{\SSS}{\maths{S}}

\newcommand{\HH}{\maths{H}}

\newcommand{\ZZ}{\maths{Z}}
\newcommand{\PP}{\maths{P}}

\newcommand{\KKK}{{\mathfrak K}}

\newcommand{\ra}{\rightarrow}
\newcommand{\bs}{\backslash}

\newcommand{\ov}[1]{{\overline #1}} 
\newcommand{\wt}[1]{{\widetilde{#1}}}

\newcommand{\ga}{\gamma}
\newcommand{\Ga}{\Gamma}

\newcommand{\cqfd}{\hfill$\Box$}

\newcommand{\card}{{\operatorname{Card}}}

\renewcommand{\Im}{{\operatorname{Im}}}
\newcommand{\Vol}{\operatorname{Vol}}

\newcommand{\arcosh}{\operatorname{argcosh}}
\newcommand{\arsinh}{\operatorname{argsinh}}

\newcommand{\bigO}{\operatorname{O}}

\renewcommand{\log}{\operatorname{ln}}

\newcommand{\CAT}{\operatorname{CAT}}

\newcommand{\hdr}{{\HH}^2_\RR}
\newcommand{\htr}{{\HH}^3_\RR}

\newcommand{\hnr}{{\HH}^n_\RR}
\newcommand{\hnc}{{\HH}^n_\CC}

\newcommand\normalout{\partial^1_{+}}

\newcommand\normalpm{\partial^1_{\pm}}

\newcounter{fig}

\usepackage{enumitem}
\setlist{nolistsep}

\title{On the hyperbolic orbital counting problem \\ 
in conjugacy classes}
\author{Jouni Parkkonen \and Fr\'ed\'eric Paulin} 
\date{\today}

\begin{document}
\bibliographystyle{../alphas}
\maketitle

\begin{abstract}
  Given a discrete group $\Ga$ of isometries of a negatively curved
  manifold $\wt M$, a non\-trivial conjugacy class $\KKK$ in $\Ga$ and
  $x_0\in\wt M$, we give asymptotic counting results, as
  $t\ra+\infty$, on the number of orbit points $\ga x_0$ at distance
  at most $t$ from $x_0$, when $\ga$ is restricted to be in $\KKK$, as well
  as related equidistribution results.  These results generalise and
  extend work of Huber on cocompact hyperbolic lattices in dimension
  $2$.  We also study the growth of given conjugacy classes in
  finitely generated groups endowed with a word metric.  
\end{abstract}

\section{Introduction}

Given an infinite discrete group of isometries $\Ga$ of a proper
metric space $X$, the {\it orbital counting problem} studies, for
fixed $x_0,y_0\in X$, the asymptotic as $t\ra +\infty$ of
$$
\card\{\ga\in\Ga\;:\; d(x_0,\ga y_0)\leq t\}\;.
$$
Initiated by Gauss in the Euclidean plane and by Huber in the real
hyperbolic plane, there is a huge corpus of works on this problem,
including the seminal results of Margulis's thesis, see for instance
\cite{Babillot02a,Oh10,Oh13} and their references for historical
remarks, as well as \cite{AthBufEskMir12, PauPolSha, Quint05,
  Sambarino13} for variations.

Given an infinite subset of the orbit $\Ga x_0$, defined in either an
algebraic, a geometric or a probabilistic way, it is interesting to
study the asymptotic growth of this subset, see for example
\cite{PetRis09,BouKonSar10} and Chapter 4 of \cite{PauPolSha} for
recent examples.  In this paper, we consider the orbit points under
the elements of a fixed nontrivial conjugacy class $\KKK$ in
$\Ga$. More precisely, we will study the asymptotic growth as
$t\ra+\infty$ of the counting function
$$
N_{\KKK,\,x_0}(t)=\card\{\ga\in\KKK\;:\; d(x_0,\ga x_0)\leq t\}\;
$$
introduced by Huber \cite{Huber56} in a special case.  Although we
will work in the framework of negative curvature in this paper, the
counting problem in (infinite) conjugacy classes is interesting even
for discrete isometry groups in (nonabelian) nilpotent or solvable Lie
groups endowed with left-invariant distances. We refer to Section
\ref{sect:countingfingengroup} for examples of computations of the
growth of $N_{\KKK,\,x_0}(t)$ when $\Ga$ is a finitely generated group
and $X$ is the set $\Ga$ endowed with a word metric.  This paper opens
a new field of research, studying which growth types (or relative
growth types) fixed conjugacy classes may have in finitely generated
groups. For word-hyperbolic groups and negatively curved manifolds,
the conjugacy classes usually have constant exponential growth rate,
as illustrated by the following result (see also Proposition
\ref{prop:growthconjclaswordhyp} and Corollary \ref{coro:encadreloxo}
for generalisations).

\btheo\label{theo:loggrowthintro} 
If $M$ is a compact negatively curved Riemannian manifold, if
$h$ is the topological entropy of its geodesic flow, if $\Ga$ is the
covering group of a universal Riemannian cover $X\ra M$, if $\KKK$ is
a nontrivial conjugacy class in $\Ga$, then
$$
\lim_{t\ra+\infty} \frac{1}{t}\log N_{\KKK,\,x_0}(t)=\frac{h}{2}\;.
$$
\etheo

In this introduction from now on, we concentrate on the case when $X$
is the real hyperbolic plane $\HH^2_\RR$, and we assume that $x_0$ is
not fixed by any nontrivial element of $\Ga$, see the main body of the
text for more general statements. Given a nontrivial element $\ga$ of
a discrete group of isometries $\Ga$ of $\HH^2_\RR$, we will denote by
$C_{\ga},\tau_{\ga},\iota_\ga$ the following objects:
\smallskip
\begin{itemize}
\item[$\bullet$] if $\ga$ is loxodromic, then $C_{\ga}$ is the translation
axis of $\ga$; with $\ell_\ga$ the translation length of $\ga$, we
define $\tau_{\ga}=(\frac{\cosh\ell_{\ga}\,-1}{2})^{1/2}$ if $\ga$
preserves the orientation and $\tau_{\ga}= (\frac{\cosh \ell_{\ga}
  \,+1}{2})^{1/2}$ otherwise; $\iota_\ga$ is $2$ if there exists an
element in $\Ga$ exchanging the two fixed points of $\ga$, and $1$
otherwise;
\item[$\bullet$]
 if $\ga$ is parabolic, then $C_{\ga}$ is a horoball centred
at the parabolic fixed point of $\ga$; we set $\tau_{\ga}= 2\sinh
\frac{d(x,\ga x)}{2}$ for any $x\in \partial C_{\ga}$; we define
$\iota_\ga$ as $2$ if there exists a nontrivial elliptic element of
$\Ga$ fixing the fixed point of $\ga$, and $1$ otherwise;
\item[$\bullet$]
 if $\ga$ is elliptic, then $C_{\ga}$ is the fixed point set
of $\ga$ in $\HH^2_\RR$; if $\ga$ is orientation-reversing, we assume
in this introduction that the stabiliser of $C_{\ga}$ is infinite;
we set $\tau_{\ga}=\sin\frac{\theta}{2}$ if $\ga$ preserves the
orientation with rotation angle $\theta$, and $\tau_{\ga}=1$
otherwise; we define $\iota_\ga=1$, unless $\ga$ preserves the
orientation with rotation angle different from $\pi$ and the
stabiliser in $\Ga$ of $C_\ga$ is dihedral, in which case
$\iota_\ga=2$.
\end{itemize}

\smallskip
We refer for instance to \cite{Roblin03} for the definition of the
critical exponent $\delta_\Ga$ of $\Ga$, the Patterson-Sullivan
measures $(\mu_{x})_{x\in \HH^2_\RR}$ of $\Ga$, the Bowen-Margulis
measure $m_{\rm BM}$ of $\Ga$, and to \cite{OhSha13,ParPau13ETDS} for
the definition of the skinning measure $\sigma_C$ of $\Ga$ associated
to a nonempty proper closed convex subset $C$ of $\HH^2_\RR$ (see also
Section \ref{sec:rappels}). We denote by $\|\mu\|$ the total mass of a
measure $\mu$ and by $\Delta_x$ the unit Dirac mass at a point $x$.

The following result says in particular that the exponential growth
rate of the orbit under a conjugacy class is $\frac{\delta_\Ga}{2}$
and that the unit tangent vectors at $x_0$ to these orbit points
equidistribute to the pulled-back Patterson-Sullivan measure.

\btheo \label{theo:intro}
Let $\Ga$ be a nonelementary finitely generated discrete group
of isometries of $\HH^2_\RR$, and let $\KKK$ be the conjugacy class of
a fixed nontrivial element $\ga_0\in\Ga$. 

\smallskip (1) As $t\ra+\infty$, we have 
$$
N_{\KKK,\,x_0}(t)\sim
\frac{\iota_{\ga_0}\,\|\mu_{x_0}\|\;\|\sigma_{C_{\ga_0}}\|}
{\delta_\Ga\,\|m_{\rm BM}\|\,{\tau_{\ga_0}}^{\delta_\Ga}}\;
e^{\frac{\delta_\Ga}{2}\;t}\;.
$$
If $\Ga$ is arithmetic or if $M$ is compact, then
the error term is $\bigO(e^{(\frac{\delta_\Ga}{2}-\kappa) t})$ for some $\kappa
>0$.  

\smallskip (2) Let $v_\ga$ be the unit tangent vector at $x_0$ to the
geodesic segment $[x_0,\ga x_0]$ for every nontrivial $\ga\in\Ga$, and
let $\pi_+:T^1_{x_0}\HH^2_\RR\ra \partial_\infty \HH^2_\RR$ be the
homeomorphism sending $v$ to the point at infinity of the geodesic ray
with initial vector $v$. For the weak-star convergence of measures on
$T^1_{x_0}\HH^2_\RR$, we have
$$
\lim_{t\ra+\infty}\; \frac{\delta_\Ga\;\|m_{\rm BM}\|
  \,{\tau_{\ga_0}}^{\delta_\Ga}}{\|\mu_{x_0}\|\;\|\sigma_{C_{\ga_0}}\|} 
\;e^{-\frac{\delta_\Ga}{2}\; t} \sum_{\ga\in\KKK,\; d(x_0,\ga x_0)\leq t}\; 
\Delta_{v_\ga} = (\pi_+^{-1})_*\mu_{x_0}\;.
$$
\etheo

When $\Ga$ is a cocompact lattice in dimension $2$ and $\ga_0$ is
loxodromic, the first claim is due to Huber \cite[Theorem B]{Huber56}
with an improved error bound in \cite{Huber98}. The following
corollary, proved in Sections \ref{sect:loxodromic} and
\ref{sect:parabolic}, is a generalisation of Huber's result for
noncompact quotients and for parabolic conjugacy classes. A version
for elliptic conjugacy classes follows from Corollary
\ref{coro:mainellip}, we leave the formulation for the reader.

\bcoro\label{coro:intro} Let $\Ga$ be a torsion-free and
orientation-preserving discrete group of isometries of $\HH^2_\RR$
such that the surface $\Ga\bs\hdr$ has finite area, with genus $g$ and
$p$ punctures.  If $\KKK$ is a conjugacy class of primitive loxodromic
elements with translation length $\ell$, then as $t\ra+\infty$,
$$ 
N_{\KKK,\,x_0}(t)\sim
\frac{\ell}{2\pi(2g+p-2)\sinh\frac\ell 2}\;e^{\frac t2}\,.
$$
If $\KKK$ is a conjugacy class of primitive parabolic elements, then
as $t\ra+\infty$,
$$
N_{\KKK,\,x_0}(t) \sim\frac{1}{2\pi(2g+p-2)}\;e^{\frac t2}\,.
$$
\ecoro

\medskip When $\Ga$ is a uniform lattice, $\KKK$ is a conjugacy class of
loxodromic elements, and $\HH^2_\RR$ is replaced by a regular tree,
the analog of Corollary \ref{coro:intro} is due to \cite{Douma11}. See
\cite{BroParPau13} for the case of any locally finite tree and more
general discrete groups of isometries.

\medskip The main tool of this paper (see Section \ref{sec:rappels})
is a counting and equidistribution result for the common
perpendiculars between locally convex subsets of simply connected
negatively curved manifolds, proved in \cite{ParPau14}. In Section
\ref{sect:counting}, we will use this tool in order to prove our
abstract main result, Theorem \ref{theo:genHuberabstrait}, on the
counting of the orbit points by the elements in a given conjugacy
class.  In Sections \ref{sect:loxodromic}, \ref{sect:parabolic} and
\ref{sect:elliptic}, we give the elementary computations concerning
the geometry of loxodromic, parabolic and elliptic isometries of a
simply connected negatively curved manifold required to apply our
abstract main result, proving as a special case the above Theorem
\ref{theo:intro}. Finally, in Section \ref{sect:countingsubgroups}, we
give some results on the counting problem of subgroups of $\Ga$ in a
given conjugacy class of subgroups.

%In the final paragraph of his introduction, Huber writes about 
%arithmetic applications via integral ternary quadratic forms 
%with reference to Fricke-Klein

\section{Counting in conjugacy classes in finitely 
generated groups}
\label{sect:countingfingengroup}

In this section, we study the growth of a given conjugacy class in a
finitely generated group endowed with a word metric, by giving three
examples. We thank E.~Breuillard, Y.~Cornulier, S.~Grigorchuk, D.~Osin
and R.~Tessera for discussions on this topic.

Let $\Ga$ be a finitely generated group, endowed with a finite
generating set $S$. For every $\ga\in\Ga$, we denote by $\|\ga\|$ the
smallest length of a word in $S\cup S^{-1}$ representing $\ga$. We
endow $\Ga$ with the left-invariant word metric $d_S$ associated to
$S$, that is, $d_S(\ga,\ga')=\|\ga^{-1}\ga'\|$ for all $\ga,\ga'
\in\Ga$. Given a conjugacy class $\KKK$ in $\Ga$, we want to study the
growth as $n\ra+\infty$ of
$$
N_\KKK(n)=N_{\KKK,\,S}(n)=\card\;\KKK\cap B_{d_S}(e,n)\;,
$$
the cardinality of the intersection of the conjugacy class $\KKK$ with
the ball of radius $n$ centered at the identity element $e$ for the
word metric $d_S$. 

Given two maps $f,g:\NN\ra\mathopen{]}0,+\infty\mathclose{[}\,$, we
write $f\asymp g$ if there exists $c\in\NN-\{0\}$ such that $g(n)\leq
f(c\,n)$ and $f(n)\leq g(c\,n)$ for every $n\in\NN$. Note that if $S'$ is
another finite generating set of $\Ga$, then $N_{\KKK,\,S'} \asymp
N_{\KKK,\,S}$.

The growth of a given conjugacy class in $\Ga$ is at most the growth
of $\Ga$, and we refer for instance to \cite{Grigorchuk14, Mann12}
and their references for information on the growth of groups.  The
growth of the trivial conjugacy class is trivial ($N_{\{e\}}(n)=1$ for
every $n\in\NN$). It would be interesting to know what are the
possible growths of given conjugacy classes, between these two
extremal bounds, and for which group all nontrivial conjugacy classes
have the same growth. We only study two examples below.

The counting problem introduced in this paper is dual to the study of
the asymptotic as $n\ra+\infty$ either of the number of translation
axes of (primitive) loxodromic elements meeting the ball of center
$x_0$ and radius $n$, in the negatively curved manifold case, or of
the number of (primitive) conjugacy classes meeting the ball of radius
$n$, in the finitely generated group case. These asymptotics have been
studied a lot, for instance by Bowen and Margulis in the manifold
case, and by Hull-Osin \cite{HulOsi13} (see also the references of
\cite{HulOsi13}) in the finitely generated group case. In particular,
Ol'shanskii \cite[Theo.~41.2]{Olshanskii91} has constructed groups
with exponential growth rate and only finitely many conjugacy classes:
at least one of them has the same growth rate as the whole group,
contrarily to the examples below.

\medskip First, let $\Ga=F(S)$ be the free group on a finite set $S$
of cardinality $|S|\geq 2$. Let $\KKK$ be the conjugacy class in $\Ga$
of a nonempty reduced and cyclically reduced word in $S\cup S^{-1}$,
denoted by $\ga_0$, of length $\ell_\KKK= \inf_{\ga\in\KKK} \|\ga\|$.
Let $m_\KKK$ be the number of cyclic conjugates of $\ga_0$ (for
instance $m_\KKK=1$ if $ \ga_0=s^\ell$ for some $s\in S\cap S^{-1}$).
We denote by $\lfloor x\rfloor$ the integral part of a real number
$x$.

\bprop\label{prop:freegroup} 
For every $n\in\NN$ with $n\geq \ell_\KKK+2$, we have
$$
N_{\KKK,\,S}(n)=m_\KKK\;(2\,|S|-2)\;
(2\,|S|-1)^{\big\lfloor\frac{n-\ell_\KKK-2}{2}\big\rfloor}\;.
$$
\eprop

In particular, $\lim_{n\ra+\infty} \frac{1}{n}\log N_{\KKK,\,S}(n)=
\frac{\log (2|S|-1)}{2}$ does not depend on the nontrivial conjugacy
class $\KKK$, and is half the exponential growth rate of $\Ga$ with
respect to the generating set $S$ (see Proposition
\ref{prop:growthconjclaswordhyp} for a generalisation).

\medskip \dem Let $k=|S|$ and $\ell=\ell_\KKK$. Every element $\ga$ in
$\KKK$ can be written uniquely as $\alpha\ga_0'\alpha^{-1}$, where
$\alpha$ is a reduced word in $S\cup S^{-1}$ and $\ga'_0$ is a cyclic
conjugate of $\ga_0$, and where the writing is reduced, that is, the
last letter of $\alpha$ is different from the inverse of the first
letter $s_1$ of $\ga_0'$ and from the last letter $s_\ell$ of
$\ga'_0$. In particular,
$$
\|\ga\|=2\,\|\alpha\|+\ell\;.
$$
Note that $s_1^{-1}\neq s_\ell$, since $\ga_0$ is cyclically
reduced. For every $m\in\NN$ with $m\geq 1$, there are
$(2k-2)(2k-1)^{m-1}$ reduced words of length at most $m$ whose last
letter is different from $s_1^{-1}$ and $s_\ell$. The result
follows. \cqfd

\medskip \rem The group $\Ga=F(S)$ acts faithfully on its Cayley graph
associated to $S$ by left multiplication, and
$N_{\KKK,\,S}(n)=\card(\KKK\cdot e\cap B(e,n))$. Proposition
\ref{prop:freegroup} gives an exact expression for this orbit count,
improving \cite[Thm.~1]{Douma11} in this special case for
$(q+1)$-regular trees with $q$ odd.

\medskip The following result says in particular that in a
torsion-free word-hyperbolic group, the nontrivial conjugacy classes
have constant exponential growth rate, equal to half the one of the
ambient group. Recall (see for instance \cite[\S 5.1]{Champetier00})
that the {\it virtual center} $Z^{\rm virt}(\Ga)$ of a nonelementary
word-hyperbolic group $\Ga$ is the finite subgroup of $\Ga$ consisting
of the elements $\ga\in\Ga$ acting by the identity on the boundary at
infinity $\partial_\infty\Ga$ of $\Ga$, or, equivalently, having
finitely many conjugates in $\Ga$, or, equivalently, whose centraliser
in $\Ga$ has finite index in $\Ga$. Note that $N_{\KKK}(n)$ is bounded
if (and only if) $\KKK$ is the conjugacy class of an element in the
virtual center.

\bprop 
\label{prop:growthconjclaswordhyp}
Let $\Ga$ be a nonelementary word-hyperbolic group, $S$ a finite
generating set of $\Ga$, and $\KKK$ the conjugacy class of an element
in $\Ga-Z^{\rm virt}(\Ga)$. Then
$$
\limsup_{n\ra+\infty} \frac{1}{n}\log N_{\KKK,\,S}(n)=
\frac{1}{2}\;\limsup_{n\ra+\infty} \frac{1}{n}\log\card\; B_{d_S}(e,n)\;.
$$
\eprop

\dem Let $\ga_0\in \Ga-Z^{\rm virt}(\Ga)$ and $\delta=
\limsup_{n\ra+\infty} \frac{1}{n}\log\card\; B_{d_S}(e,n)$. Let
$C_{\ga_0}$ be a quasi-translation axis of $\ga_0$ if $\ga_0$ has
infinite order, and let $C_{\ga_0}$ be the set of quasi-fixed points
of $\ga_0$ otherwise. Note that $C_{\ga_0}$ is quasi-convex, that
$Z_\Ga(\ga_0)$ preserves $C_{\ga_0}$, and that $C_{\ga_0}$ is at
bounded distance from $Z_\Ga(\ga_0)$ in $\Ga$. In particular,
$\Ga_0=Z_\Ga(\ga_0)$ is a quasi-convex-cocompact subgroup with
infinite index in the nonelementary word hyperbolic group $\Ga$. 

It is well-known that the exponential growth rate of $\Ga/\Ga_0$ is
then equal to the exponential growth rate $\delta$ of $\Ga$. Indeed,
the limit set $\Lambda\Ga_0$ of $\Ga_0$ is then a proper subset of
$\partial_\infty\Ga$ and $\Ga_0$ acts properly discontinuously on
$\Ga\cup (\partial_\infty\Ga-\Lambda\Ga_0)$. Let $\xi
\in \partial_\infty \Ga- \Lambda\Ga_0$.  If $U$ is a small enough
neighbourhood of $\xi$ in $\Ga\cup\partial_\infty\Ga$, then there
exists $N\in\NN$ such that $U$ meets at most $N$ of its images by the
elements of $\Ga_0$, and for every $x\in U\cap \Ga$, if
$p:\Ga\ra\Ga/\Ga_0$ is the canonical projection, then
$|d(x,e)-d(p(x),p(e))|$ is uniformly bounded. It is well-known (see
for instance the proof of \cite[Coro.~1]{Roblin02}) that the
(sectorial) exponential growth rate $\limsup_{n\ra+\infty}
\frac{1}{n}\log\card\; \big(U\cap B_{d_S}(e,n)\big)$ of $\Ga$ in $U$
is equal to $\delta$. This proves the above claim.

Up to a bounded additive constant, the distance between $e$ and
$\ga^{-1}\ga_0\ga$ is equal to twice the distance from $\ga$ to
$C_{\ga_0}$, by hyperbolicity. Hence the exponential growth rate of
$\KKK$ is half the exponential growth rate of $\Ga/Z_\Ga(\ga_0)$, that
is $\delta/2$.  \cqfd

\bigskip Now, let $A$ be a free abelian group of rank $2k$, let
$\langle\cdot,\cdot\rangle$ be an integral symplectic form on $A$, and
let $\Ga$ be the  associated Heisenberg group, that is, the group
with underlying set $A\times\ZZ$ and group law
$$
(a,z)(a',z')=(a+a',z+z'+\langle a,a'\rangle)\;.
$$
Note that $\Ga$ is finitely generated, and we have an exact sequence
of groups
$$
0\longrightarrow \ZZ\stackrel{i}{\longrightarrow}\Ga
\stackrel{\pi}{\longrightarrow}A\longrightarrow 0\;,
$$
where $i:z\mapsto (0,z)$ has image the center of $\Ga$ and
$\pi:(a,z)\mapsto a$. Let $\KKK$ be a nontrivial conjugacy class (that
is, the conjugacy class of a noncentral element) in $\Ga$.

\bprop We have
$$
N_{\KKK}(n)\asymp n^2\;.
$$
\eprop

In particular, the growth of any nontrivial conjugacy class in the
Heisenberg group $\Ga$ is quadratic. Note that $\card\;B_\Ga(e,n)
\asymp n^{2k+2}$ and that the number of (primitive or not) conjugacy
classes meeting the ball of radius $n$ is $\asymp n^2\ln n$ if $k=1$,
see \cite[Ex.~2.4]{GubSap10}.

\medskip \dem Let $\ga_0=(a_0,z_0)$ be a noncentral element in $\Ga$,
so that $a_0\neq 0$, and let $\|\ga_0\|$ be its distance to the
identity element $e$ for a given word metric on $\Ga$.

Since $\pi:\Ga\ra A$ is the abelianisation map, whose kernel is the
center $Z$ of $\Ga$, we have $\pi(\KKK)=\{\pi(\ga_0)\}$ and
$$
\KKK\subset \pi^{-1}(\{\pi(\ga_0)\})=Z\,\ga_0\;.
$$

Since $\langle\cdot,\cdot\rangle$ is nondegenerate and $a_0\neq 0$,
there exists $b_0\in A$ such that $n_0=2\langle a_0,b_0\rangle\neq
0$. For every $(a,z)\in\Ga$, since $(a,z)^{-1}=(-a,-z)$, it is easy to
compute that
$$
(a,z)(a_0,z_0)(a,z)^{-1}=(a_0,z_0+2\langle a,a_0\rangle)\;.
$$
Hence, with $Z^{n_0}=\{(0,n_0 n)\;:\;n\in\ZZ\}$, which is a finite
index subgroup of $Z$, we have
$$ 
Z^{n_0}\,\ga_0\subset \KKK\;.
$$
We have
$$
\card\;\KKK\cap B(e,n) \leq
\card\,\big(Z\cap B(e,n)\,\ga_0^{-1}\big) \leq
\card\;Z\cap B(e,n+\|\ga_0\|)\;,
$$
and similarly, $\card\;\KKK\cap B(e,n) \geq\card\;Z^{n_0}\cap
B(e,n-\|\ga_0\|)$.  We hence only have to prove that for every finite
index subgroup $Z'$ of $Z$, we have $\card\;Z\cap B(e,n)\asymp
n^2$. This is well known (see for instance \cite[VII.21]{Harpe00} when
$A$ has rank $2$): for instance, we have $[(a_0,0)^p,(b_0,0)^q]
=(0,pq\, \langle a_0,b_0\rangle)$ for all $p,q\in\ZZ$, and the
distance to $e$ of the commutator on the left hand side of this
equality is at most $c(p+q)$, for some constant $c>0$.  \cqfd

\section{Counting common perpendicular arcs}
\label{sec:rappels}

In this section, we briefly review a simplified version of the
geometric counting and equidistribution result proved in
\cite{ParPau14}, which is the main tool in this paper (see also
\cite{ParPauRev} for related references, \cite{ParPauTou} for
arithmetic applications in real hyperbolic spaces and
\cite{ParPauHeis} for the case of locally symmetric spaces).
We refer to \cite{BriHae99} for the background definitions and
properties concerning the isometries of $\CAT(-1)$ spaces.

Let $\wt M$ be a complete simply connected Riemannian manifold with
(dimension at least $2$ and) pinched negative sectional curvature
$-b^2\le K\le -1$, let $x_0\in\wt M$, and let $T^1\wt M$ be the unit
tangent bundle of $\wt M$. Let $\Ga$ be a nonelementary discrete group
of isometries of $\wt M$ and let $M=\Ga\bs\wt M$ and $T^1M= \Ga\bs
T^1\wt M$ be the quotient orbifolds.

We denote by $\partial_\infty \wt M$ the boundary at infinity
of $\wt M$, by $\Lambda \Ga$ the limit set of $\Ga$ and by
$(\xi,x,y)\mapsto \beta_\xi(x,y)$ the Busemann cocycle on
$\partial_\infty \wt M\times \wt M\times \wt M$ defined by
$$
(\xi, x,y)\mapsto \beta_{\xi}(x,y)=
\lim_{t\to+\infty}d(\rho_t,x)-d(\rho_t,y)\;,
$$
where $\rho:t\mapsto \rho_t$ is any geodesic ray with point at
infinity $\xi$ and $d$ is the Riemannian distance.  

For every $v\in T^1\wt M$, let $\pi(v)\in \wt M$ be its origin, and
let $v_-, v_+$ be the points at infinity of the geodesic line in $\wt
M$ whose tangent vector at time $t=0$ is $v$.  We denote by
$\pi_+:T^1_{x_0}\wt M\ra \partial_\infty \wt M$ the homeomorphism
$v\mapsto v_+$.

\medskip Let $D^-$ and $D^+$ be nonempty proper closed convex subsets
in $\wt M$, with stabilisers $\Ga_{D^-}$ and $\Ga_{D^+}$ in $\Ga$,
such that the families $(\ga D^-)_{\ga\in\Ga/\Ga_{D^-}}$ and $(\ga
D^+)_{\ga\in\Ga/\Ga_{D^+}}$ are locally finite in $\wt M$. We denote
by $\normalpm D^\mp$ the {\it outer/inner unit normal bundle} of
$\partial D^\mp$, that is, the set of $v\in T^1\wt M$ such that
$\pi(v)\in \partial D^\mp$, $v_\pm \in \partial_\infty \wt
M-\partial_\infty D^\mp$ and the closest point projection on $D^\mp$
of $v_\pm $ is $\pi(v)$. For every $\ga,\ga'$ in $\Ga$ such that $\ga
D^-$ and $\ga' D^+$ have a common perpendicular (that is, if the
closures $\overline{\ga D^-}$ and $\overline{\ga' D^+}$ in $\wt
M\cup\partial_\infty \wt M$ are disjoint), we denote by
$\alpha_{\ga,\,\ga'}$ this common perpendicular (starting from $\ga
D^-$ at time $t=0$), by $\ell(\alpha_{\ga,\,\ga'})$ its length, and by
$v^-_{\ga,\,\ga'} \in \ga \normalout D^-$ its initial tangent vector.
The {\em multiplicity} of $\alpha_{\ga,\,\ga'}$ is
$$
m_{\ga,\,\ga'}=
\frac 1{\card(\ga\Ga_{D^-}\ga^{-1}\cap\ga'\Ga_{D^+}{\ga'}^{-1})}\,,
$$
which equals $1$ when $\Ga$ acts freely on $T^1\wt M$ (for instance
when $\Ga$ is torsion-free).  Let
$$
\N_{D^-,\,D^+}(s)=\sum_{\substack{
(\ga,\,\ga')\in \,\Ga\bs((\Ga/\Ga_{D^-})\times (\Ga/\Ga_{D^+}))\\
\phantom{\big|}\overline{\ga D^-}\,\cap \,\overline{\ga' D^+}\,
=\emptyset,\; \ell(\alpha_{\ga,\, \ga'})\leq s}} m_{\ga,\,\ga'}=
\sum_{\substack{[\ga]\in\, \Ga_{D^-}\bs\Ga/\Ga_{D^+}\\
\phantom{\big|}\overline{D^-}\,\cap \,\overline{\ga D^+}\,=
\emptyset,\; \ell(\alpha_{e,\, \ga})\leq s}} m_{e,\,\ga}
\;,
$$
where $\Ga$ acts diagonally on $(\Ga/\Ga_{D^-})\times
(\Ga/\Ga_{D^+})$. When $\Ga$ is torsion-free, $\N_{D^-,\,D^+}(s)$ is
the number of the common perpendiculars of length at most $s$ between
the images of $D^-$ and $D^+$ in $M$, with multiplicities coming from
the fact that $\Ga_{D^\pm}\bs D^\pm$ is not assumed to be embedded in
$M$. We refer to \cite[\S 4]{ParPau14} for the use of
H\"older-continuous potentials on $T^1\wt M$ to modify this counting
function by adding weights.

Recall the following notions (see for instance \cite{Roblin03}). The
{\em critical exponent} of $\Ga$ is
$$
\delta_{\Ga}=\limsup_{N\to+\infty}\frac 1N \ln\card\{\ga\in\Ga:
d(x_{0},\ga x_{0})\leq N\}\,,
$$
which is positive, finite, independent of $x_{0}$ (and equal to the
topological entropy $h$ if $\Ga$ is cocompact and torsion-free). Let
$(\mu_{x})_{x\in \wt M}$ be a {\em Patterson-Sullivan density} for
$\Ga$, that is, a family $(\mu_{x})_{x\in \wt M}$ of nonzero finite
measures on $\partial_{\infty}\wt M$ whose support is $\Lambda\Ga$,
such that $\ga_*\mu_x=\mu_{\ga x}$ and
$$
\frac{d\mu_{x}}{d\mu_{y}}(\xi)=e^{-\delta_{\Ga}\beta_{\xi}(x,\,y)}
$$
for all $\ga\in\Ga$, $x,y\in\wt M$ and $\xi\in\partial_{\infty}\wt M$.
The {\em Bowen-Margulis measure} $\wt m_{\rm BM}$ for $\Ga$ on $T^1\wt
M$ is defined, using Hopf's parametrisation $v\mapsto
(v_-,v_+,\beta_{v_+}(x_0,\pi(v))\,)$ of $T^1\wt M$, by
$$
d\wt m_{\rm BM}(v)=e^{-\delta_{\Ga}(\beta_{v_{-}}(\pi(v),\,x_{0})+
\beta_{v_{+}}(\pi(v),\,x_{0}))}\;
d\mu_{x_{0}}(v_{-})\,d\mu_{x_{0}}(v_{+})\,dt\,.  
$$
The measure $\wt m_{\rm BM}$ is nonzero and independent of $x_{0}$.
It is invariant under the geodesic flow, the antipodal map
$\iota:v\mapsto -v$ and the action of $\Ga$, and thus defines a
nonzero measure $m_{\rm BM}$ on $T^1M=\Ga\bs T^1\wt M$, called the
{\em Bowen-Margulis measure} on $M$, which is invariant under the
geodesic flow of $M$ and the antipodal map. If $\wt M$ is symmetric
and if $\Ga$ is geometrically finite (for instance if $\wt M=\hdr$ and
$\Ga$ is finitely generated), then $m_{\rm BM}$ is finite. See for
instance \cite{DalOtaPei00} for many more examples. If $m_{\rm BM}$ is
finite, then $m_{\rm BM}$ is mixing under the geodesic flow if $\wt M$
is symmetric or if $\Lambda\Ga$ is not totally disconnected (hence if
$M$ is compact), see for instance \cite{Babillot02b,DalBo99}.

\medskip Using the endpoint homeomorphisms $v\mapsto v_\pm$ from
$\normalpm {D^\mp}$ to $\partial_{\infty}\wt M-\partial_{\infty}D^\mp$,
the {\em skinning measure}
$\wt\sigma_{D^\mp}$ of $\Ga$ on $\normalpm{D^\mp}$ is defined by
$$
d\wt\sigma_{D^\mp}(v) =  
e^{-\delta\,\beta_{v_{\pm}}(\pi(v),\,x_{0})}\,d\mu_{x_{0}}(v_{\pm})\,,
$$
see \cite[\S 1.2]{OhSha13} when $D^\mp$ is a horoball or a totally
geodesic subspace in $\wt M$ and \cite{ParPau13ETDS},
\cite{ParPau14} for the general case of convex subsets in variable
curvature and with a potential. 

The measure $\wt\sigma_{D^\mp}$ is independent of $x_{0}\in\wt M$, it
is nonzero if $\Lambda\Ga$ is not contained in $\partial_{\infty}
D^\mp$, and satisfies $\wt\sigma_{\ga D^\mp} =\ga_{*}\wt
\sigma_{D^\mp}$ for every $\ga\in\Ga$. Since the family $(\ga
D^\mp)_{\ga\in\Ga/\Ga_{D^\mp}}$ is locally finite in $\wt M$, the
measure $\sum_{\ga \in \Ga/\Ga_{D^\mp}} \;\ga_*\wt\sigma_{D^\mp}$ is a
well defined $\Ga$-invariant locally finite (nonnegative Borel)
measure on $T^1\wt M$, hence it induces a locally finite measure
$\sigma_{D^\mp}$ on $T^1M$, called the {\em skinning measure} of
$D^\mp$ in $T^1M$. If $\Ga_{D^\mp}\bs\partial D^\mp$ is compact, then
$\sigma_{D^\mp}$ is finite. We refer to \cite[\S 5]{OhShaCircles} and
\cite[Theo.~9]{ParPau13ETDS} for finiteness criteria of the skinning
measure $\sigma_{D^\mp}$.

\medskip The following result on the asymptotic behaviour of the
counting function $\N_{D^-,\,D^+}$ is a special case of more general
results \cite[Coro.~20, 21, Theo.~28]{ParPau14}. We refer to
\cite{ParPauRev} for a survey of the particular cases known before
\cite{ParPau14}, due to Huber, Margulis, Herrmann, Cosentino, Roblin,
Eskin-McMullen, Oh-Shah, Martin-McKee-Wambach, Pollicott, and the
authors for instance.

\btheo\label{theo:mainequicountdown} Let $\Ga,D^-,D^+$ be as
above. Assume that the measures $m_{\rm BM},\sigma_{D^-},\sigma_{D^+}$
are nonzero and finite, and that $m_{\rm BM}$ is mixing for the
geodesic flow of $T^1M$. Then
$$
\N_{D^-,\,D^+}(s)\;\sim\;
\frac{\|\sigma_{D^-}\|\;\|\sigma_{D^+}\|}{\delta_\Ga\;\|m_{\rm
    BM}\|}\; e^{\delta_\Ga \,s}\;,
$$
as $s\ra+\infty$. If $\Ga$ is arithmetic or if $M$ is compact, then
the error term is $\bigO(e^{(\delta_\Ga-\kappa) s})$ for some $\kappa
>0$.  Furthermore, the initial vectors of the common perpendiculars
equidistribute in the outer unit normal bundle of $D^-$: 
\begin{equation}\label{eq:equidistribdown}
\lim_{s\ra+\infty}\; 
\frac{\delta_\Ga\;\|m_{\rm BM}\|}{\|\sigma_{D^-}\|\;\|\sigma_{D^+}\|}
\;e^{-\delta_\Ga\, s}\;
\sum_{\substack{[\ga]\in\, \Ga_{D^-}\bs\Ga/\Ga_{D^+}\\
\phantom{\big|}\overline{D^-}\,\cap \,\overline{\ga D^+}\,=
\emptyset,\; \ell(\alpha_{e,\, \ga})\leq s}} m_{e,\,\ga} \;
\Delta_{v^-_{e,\,\ga}}\;=\; \frac{\wt\sigma_{D^-}}{\|\sigma_{D^-}\|}\;
\end{equation}
for the weak-star convergence of measures on the locally compact space
$T^1\wt M$. 
\cqfd
\etheo

\section{Counting in conjugacy classes}
\label{sect:counting}

Let $\wt M,x_0,\Ga$ be as in the beginning of Section
\ref{sec:rappels}.
For any nontrivial element $\ga$ in $\Ga$, let  
$C_\ga$ be

$\bullet$~ the translation axis of $\ga$ if $\ga$ is loxodromic,

$\bullet$~ the fixed point set of  $\ga$ if $\ga$ is elliptic,

$\bullet$~ a horoball centered at the fixed point of $\ga$ if $\ga$ is
parabolic,

\noindent  which is a nonempty proper closed convex subset of $\wt M$. We assume (this condition is automatically satisfied unless $\ga'$ is parabolic)
that $\ga C_{\ga'}=C_{\ga\ga'\ga^{-1}}$ for all $\ga\in\Ga$ and
$\ga'\in \Ga- \{e\}$.  

By the equivariance properties of the skinning measures, the total
mass of $\sigma_{C_\ga}$ depends only on the conjugacy class $\KKK$ of
$\ga$, and will be denoted by $\|\sigma_\KKK\|$. This quantity, called
the {\it skinning measure} of $\KKK$, is positive and finite for
instance when $\ga$ is loxodromic, since $\Lambda\Ga$ contains at
least $3$ points and the image of $C_\ga$ in $M$ is compact. 
In Sections \ref{sect:parabolic} and \ref{sect:elliptic}, we will
give  other
classes of examples of conjugacy classes $\KKK$ with positive and
finite skinning measure $\|\sigma_\KKK\|$, and prove in particular
that this is always true if $\wt M=\HH^2_\RR$ except possibly when
$\ga$ is elliptic and orientation-reversing.

We define 
$$
m_\ga= \frac{1}{\card(\Ga_{x_0}\cap \Ga_{C_\ga})}\,,
$$ 
which is a natural complexity of $\ga$, independent on the choice of
$C_\ga$ when $\ga$ is parabolic, and equals $1$ if the stabiliser of
$x_0$ in $\Ga$ is trivial (for instance if $\Ga$ is torsion-free).
Note that for every $\alpha\in\Ga$, the real number
$m_{\alpha\ga\alpha^{-1}}$ depends only on the double coset of
$\alpha$ in $\Ga_{x_0}\bs\Ga/\Ga_{C_{\ga}}$.

The centraliser $Z_\Ga(\ga)$ of $\ga$ in $\Ga$ is contained in the
stabiliser of $C_\ga$ in $\Ga$. The index $[\Ga_{C_\ga}:Z_\Ga(\ga)]$
depends only on the conjugacy class $\KKK$ of $\ga$; it will be
denoted by $i_\KKK$ and called the {\it index} of $\KKK$. We assume in
what follows that $i_\KKK$ is finite, which is in particular the case
if $\ga$ is loxodromic (the stabiliser of its translation axis $C_\ga$
is virtually cyclic). In Sections \ref{sect:parabolic} and
\ref{sect:elliptic}, we will give other classes of examples of
conjugacy classes $\KKK$ with finite index $i_\KKK$, and prove in
particular that this is always true if $\wt M=\HH^2_\RR$.

\medskip
We define the counting function 
$$
N_{\KKK,\,x_0}(t)=\sum_{\alpha\in\KKK,\; d(x_0,\,\alpha x_0)\leq t} m_\alpha\,.
$$
When the stabiliser of $x_0$ in $\Ga$ is trivial, we
recover the definition of the Introduction.

\medskip Let $\psi:[0,+\infty\mathclose{[}\ra[0,+\infty\mathclose{[}$
be an eventually nondecreasing map such that $\lim_{t\ra+\infty}
\psi(t) = +\infty$. We will say that a nontrivial element
$\ga_0\in\Ga$ is {\it $\psi$-equitranslating} if for every $x\in \wt
M$ at distance big enough from $C_{\ga_0}$, we have
$$
d(x,C_{\ga_0})=\psi(d(x,\ga_0 x))\,.
$$
Note that this condition depends only on the conjugacy class of
$\ga_0$.  When $\ga_0$ is parabolic, up to replacing $\psi$ by
$\psi+c$ for some constant $c\in\RR$, this condition does not depend
on the choice of the horoball $C_{\ga_0}$. In Sections
\ref{sect:loxodromic}, \ref{sect:parabolic} and \ref{sect:elliptic},
we will give several classes of examples of equitranslating
isometries, and prove in particular that every nontrivial isometry of
$\HH^2_\RR$ is equitranslating.

The following theorem is the main abstract result of this paper.

\btheo\label{theo:genHuberabstrait} Assume that the Bowen-Margulis
measure of $\Ga$ is finite and mixing for the geodesic flow on $T^1M$.
Let $\KKK$ be a conjugacy class of $\psi$-equitranslating elements of
$\Ga$ with finite index $i_\KKK$ and positive and finite skinning
measure $\|\sigma_\KKK\|$.  Then, as $t\ra+\infty$,
$$
N_{\KKK,\,x_0}(t)\sim 
\frac{i_{\KKK}\,\|\mu_{x_0}\|\,\|\sigma_{\KKK}\|}
{\delta_\Ga\,\|m_{\rm BM}\|}
\,e^{\delta_\Ga\, \psi(t)}\,.
$$
If $\Ga$ is arithmetic or if $M$ is compact, then the error term is
$\bigO(e^{(\delta_\Ga-\kappa) \psi(t)})$ for some $\kappa >0$.
Furthermore, if $v_\alpha$ is the unit tangent vector at $x_0$ to the
geodesic segment $[x_0,\alpha x_0]$ for every $\alpha\in \Ga-
\Ga_{x_0}$, for the weak-star convergence of measures on $T^1\wt M$,
we have
$$
\lim_{t\ra+\infty}\; 
\frac{\delta_\Ga\;\|m_{\rm BM}\|}{i_{\KKK}\,\|\sigma_{\KKK}\|}
\;e^{-\delta_\Ga\, \psi(t)}\;
\sum_{\alpha\in\KKK,\; 0<d(x_0,\,\alpha x_0)\leq t} 
m_\alpha\,\Delta_{v_\alpha} \;=\; (\pi_+^{-1})_*\mu_{x_0}\;.
$$
\etheo

\dem 
Let $\ga_0$ be a $\psi$-equitranslating element of $\Ga-\{e\}$
and let $\KKK=\{\ga\ga_0\ga^{-1}\;:\;\ga\in\Ga\}$ be its conjugacy
class.  Since $\wt\sigma_{\{x_0\}} =(\pi_+^{-1})_*\mu_{x_0}$ (see
\cite[\S 3]{ParPau13ETDS}), we have
$$
\|\sigma_{\{x_0\}}\|= \frac{\|\mu_{x_0}\|}{|\Ga_{x_0}|}\,.
$$
In particular, both $\|\sigma_{\{x_0\}}\|$ and $\|\sigma_{\C_{\ga_0}}
\| = \|\sigma_\KKK\|$, are positive and finite. Hence, since $\psi$ is
eventually nondecreasing, by the definition of the counting
function $\N_{D^-,\,D^+}$ for $D^-=\{x_0\}$ and $D^+=C_{\ga_0}$, and
by Theorem \ref{theo:mainequicountdown}, we have, as $t\ra+\infty$,
\begin{align*}
\sum_{\alpha\in\KKK,\; 0<d(x_0,\,\alpha x_0)\leq t} m_\alpha&\sim
\sum_{\alpha\in\KKK,\; 0<d(x_0,\,C_{\alpha})\leq \psi(t)} m_\alpha = 
\sum_{\ga\in\Ga/Z_\Ga(\ga_0),\; 0<d(x_0,\,\ga C_{\ga_0})\leq \psi(t)} 
m_{\ga\ga_0\ga^{-1}}\\ & = |\Ga_{x_0}|\,i_\KKK\;
\sum_{\ga\in\Ga_{x_0}\bs\Ga/\Ga_{C_{\ga_0}},\; 0<d(x_0,\,\ga C_{\ga_0})\leq \psi(t)} 
m_{\ga\ga_0\ga^{-1}}\\  &
=|\Ga_{x_0}|\,i_\KKK\;\N_{\{x_0\},\,C_{\ga_0}}(\psi(t))
\\  &\sim
|\Ga_{x_0}|\,i_\KKK\;\frac{\|\sigma_{\{x_0\}}\|\;\|\sigma_{C_{\ga_0}}\|}
{\delta_\Ga\;\|m_{\rm BM}\|}\; e^{\delta_\Ga \,\psi(t)}\;.
\end{align*}
The first claim of Theorem \ref{theo:genHuberabstrait}
follows.  

\medskip\noindent
\begin{minipage}{8cm} ~~~ For every $\alpha\in\KKK$, let $p_\alpha$
  be the closest point to $x_0$ on $C_\alpha$. Then $\alpha p_\alpha$
  is the closest point to $\alpha x_0$ on $C_\alpha$.
\end{minipage}
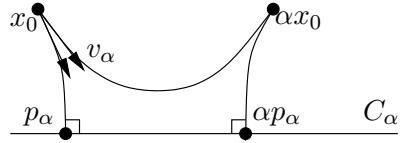
\begin{minipage}{6.9cm}
\begin{center}
\begin{picture}(0,0)%
\includegraphics{fig_initvect.pstex}%
\end{picture}%
\setlength{\unitlength}{3812sp}%
\begingroup\makeatletter\ifx\SetFigFont\undefined%
\gdef\SetFigFont#1#2#3#4#5{%
  \reset@font\fontsize{#1}{#2pt}%
  \fontfamily{#3}\fontseries{#4}\fontshape{#5}%
  \selectfont}%
\fi\endgroup%
\begin{picture}(2551,941)(1512,-1378)
\put(1527,-623){\makebox(0,0)[lb]{\smash{{\SetFigFont{11}{13.2}{\rmdefault}{\mddefault}{\updefault}{\color[rgb]{0,0,0}$x_0$}%
}}}}
\put(3252,-608){\makebox(0,0)[lb]{\smash{{\SetFigFont{11}{13.2}{\rmdefault}{\mddefault}{\updefault}{\color[rgb]{0,0,0}$\alpha x_0$}%
}}}}
\put(2026,-826){\makebox(0,0)[lb]{\smash{{\SetFigFont{11}{13.2}{\rmdefault}{\mddefault}{\updefault}{\color[rgb]{0,0,0}$v_\alpha$}%
}}}}
\put(1621,-1231){\makebox(0,0)[lb]{\smash{{\SetFigFont{11}{13.2}{\rmdefault}{\mddefault}{\updefault}{\color[rgb]{0,0,0}$p_\alpha$}%
}}}}
\put(3106,-1231){\makebox(0,0)[lb]{\smash{{\SetFigFont{11}{13.2}{\rmdefault}{\mddefault}{\updefault}{\color[rgb]{0,0,0}$\alpha p_\alpha$}%
}}}}
\put(3826,-1231){\makebox(0,0)[lb]{\smash{{\SetFigFont{11}{13.2}{\rmdefault}{\mddefault}{\updefault}{\color[rgb]{0,0,0}$C_\alpha$}%
}}}}
\end{picture}%

\end{center}
\end{minipage}

\medskip Since $\lim_{t\ra+\infty} \psi(t) = +\infty$, when $d(x_0,
\alpha x_0)$ is large enough, the distance $d(x_0,C_{\alpha})$ becomes
large.  Hence the initial tangent vector $v_{\alpha}$ to the geodesic
segment $[x_0, \alpha x_0]$ becomes arbitrarily close  to
the initial tangent vector to the geodesic segment  $[x_0,p_\alpha]$,
uniformly on $\alpha\in\KKK$ such that $d(x_0,C_{\alpha})$ is large enough, and
independently on $d(p_\alpha, \alpha p_\alpha)$ which could be $0$.
Hence, using again and similarly Theorem \ref{theo:mainequicountdown}
with $D^-=\{x_0\}$ and $D^+=C_{\ga_0}$, we have, as $t\ra+\infty$,
\begin{align*}
&\frac{\delta_\Ga\;\|m_{\rm BM}\|}{i_{\KKK}\,\|\sigma_{\KKK}\|}
\;e^{-\delta_\Ga\, \psi(t)}
\sum_{\alpha\in\KKK,\; 0<d(x_0,\,\alpha x_0)\leq t} m_\alpha\;\Delta_{v_\alpha}
\\ \sim\;\; & 
\frac{\delta_\Ga\;\|m_{\rm BM}\|}{i_{\KKK}\,\|\sigma_{\KKK}\|}
\;e^{-\delta_\Ga\, \psi(t)}
\sum_{\ga\in\Ga/Z_\Ga(\ga_0),\; 0<d(x_0,\,\ga C_{\ga_0})\leq \psi(t)} 
m_{\ga\ga_0\ga^{-1}}\;\Delta_{v_{\ga\ga_0\ga^{-1}}}
\\ \sim\;\; & 
\frac{\delta_\Ga\;\|m_{\rm BM}\|}{\|\sigma_{D^+}\|}
\;e^{-\delta_\Ga\, \psi(t)}
\sum_{\ga\in\Ga/\Ga_{D^+},\; 0<d(x_0,\,\ga D^+)\leq \psi(t)} 
m_{e,\,\ga}\;\Delta_{v^-_{e,\,\ga}}
\\ \stackrel{*}{\rightharpoonup}\;\;&
\wt\sigma_{\{x_0\}}=(\pi_+^{-1})_*\mu_{x_0}\;.
\end{align*}
This proves the second claim of Theorem \ref{theo:genHuberabstrait}.  
\cqfd

\section{The geometry of loxodromic isometries}
\label{sect:loxodromic}

In this section, we fix a loxodromic isometry $\ga$ of a complete
$\CAT(-1)$ geodesic metric space $X$. Let $\ell= \ell_{\ga}=
\inf_{x\in X} d(x,\ga x)>0$ be its translation length and let
$$
C_\ga=\{x\in X\;:\;d(x,\ga x)=\ell\}
$$ 
be its translation axis.

If $X=\hdr$, if $\ga$ is orientation-preserving, and if $x\in\hdr$ is
at a distance $s$ from the translation axis of $\ga$, then 
\begin{equation}\label{eq:planetranslation}
d(x,\ga x)=2\arsinh(\cosh s\,\sinh\frac{\ell} 2)\,.
\end{equation}
Indeed, after a conjugation by a suitable isometry, we may assume that
the translation axis of $\ga$ is the geodesic line with endpoints $0$
and $\infty$ in the upper halfplane model of $\hdr$, that $\ga
z=e^{\ell} z$ for all $z\in\CC$ with $\Im\; z>0$, and that $x$ is on
the geodesic ray starting from $i$ and ending at $1$.  Using the angle
of parallelism formula \cite[Thm.~7.9.1]{Beardon83}, we have $x=(\tanh
s,\frac 1{\cosh s})$, which gives $\ga x=e^{\ell}(\tanh s,\frac
1{\cosh s})$. From this, Equation \eqref{eq:planetranslation} follows
using the hyperbolic distance formula \cite[Thm.~7.2.1
(iii)]{Beardon83}.

In the other extreme, if $X$ is a tree and if $x\in X$ is at a
distance $s$ from the translation axis of $\ga$, then
$d(x,\ga x)=\ell+2s$.
The general situation lies between these two cases.

\blemm\label{lem:generaltranslation} If $x\in X$ is at distance
$s$ from the translation axis of $\ga$, then
$$
2\arsinh(\cosh s\sinh\frac{\ell}  2)\le d(x,\ga x)
\le 2s+\ell\,.
$$
\elemm

Note that as $s\to+\infty$, the lower bound is equal to
$2s+2\log(\sinh\frac{\ell }2) +\bigO(e^{-2s})$, hence the difference of
the upper and lower bounds is bounded by a constant that only depends
on $\ell$.

\medskip \dem The upper bound follows from the triangle inequality.
Let us prove the lower bound. Let $p$ and $q=\ga p$ be the closest
points on $C_\ga$ to respectively $x$ and $\ga x$.  Let $Q$ be the
quadrilateral in $X$ with vertices $x$, $p$, $q$ and $\ga x$.  

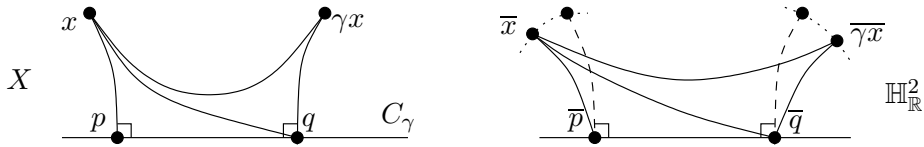
\begin{center}
\begin{picture}(0,0)%
\includegraphics{fig_quadrangl.pstex}%
\end{picture}%
\setlength{\unitlength}{3812sp}%
\begingroup\makeatletter\ifx\SetFigFont\undefined%
\gdef\SetFigFont#1#2#3#4#5{%
  \reset@font\fontsize{#1}{#2pt}%
  \fontfamily{#3}\fontseries{#4}\fontshape{#5}%
  \selectfont}%
\fi\endgroup%
\begin{picture}(5745,948)(1156,-1378)
\put(6661,-691){\makebox(0,0)[lb]{\smash{{\SetFigFont{11}{13.2}{\rmdefault}{\mddefault}{\updefault}{\color[rgb]{0,0,0}$\overline{\ga x}$}%
}}}}
\put(4827,-1248){\makebox(0,0)[lb]{\smash{{\SetFigFont{11}{13.2}{\rmdefault}{\mddefault}{\updefault}{\color[rgb]{0,0,0}$\overline{p}$}%
}}}}
\put(6886,-1096){\makebox(0,0)[lb]{\smash{{\SetFigFont{11}{13.2}{\rmdefault}{\mddefault}{\updefault}{\color[rgb]{0,0,0}$\HH^2_\RR$}%
}}}}
\put(1171,-1006){\makebox(0,0)[lb]{\smash{{\SetFigFont{11}{13.2}{\rmdefault}{\mddefault}{\updefault}{\color[rgb]{0,0,0}$X$}%
}}}}
\put(4384,-646){\makebox(0,0)[lb]{\smash{{\SetFigFont{11}{13.2}{\rmdefault}{\mddefault}{\updefault}{\color[rgb]{0,0,0}$\overline{x}$}%
}}}}
\put(6260,-1257){\makebox(0,0)[lb]{\smash{{\SetFigFont{11}{13.2}{\rmdefault}{\mddefault}{\updefault}{\color[rgb]{0,0,0}$\overline{q}$}%
}}}}
\put(1722,-1248){\makebox(0,0)[lb]{\smash{{\SetFigFont{11}{13.2}{\rmdefault}{\mddefault}{\updefault}{\color[rgb]{0,0,0}$p$}%
}}}}
\put(1527,-623){\makebox(0,0)[lb]{\smash{{\SetFigFont{11}{13.2}{\rmdefault}{\mddefault}{\updefault}{\color[rgb]{0,0,0}$x$}%
}}}}
\put(3094,-1245){\makebox(0,0)[lb]{\smash{{\SetFigFont{11}{13.2}{\rmdefault}{\mddefault}{\updefault}{\color[rgb]{0,0,0}$q$}%
}}}}
\put(3286,-601){\makebox(0,0)[lb]{\smash{{\SetFigFont{11}{13.2}{\rmdefault}{\mddefault}{\updefault}{\color[rgb]{0,0,0}$\ga x$}%
}}}}
\put(3608,-1227){\makebox(0,0)[lb]{\smash{{\SetFigFont{11}{13.2}{\rmdefault}{\mddefault}{\updefault}{\color[rgb]{0,0,0}$C_\ga$}%
}}}}
\end{picture}%

\end{center}

\noindent Let $\overline{Q}$ be the quadrilateral in $\hdr$ with
vertices $\overline{x}$, $\overline{p}$, $\overline{q}$ and
$\overline{\ga x}$, obtained by gluing together, along the geodesic
segment $[\overline{x}, \overline{q}]$, the comparison triangles of
the two triangles in $X$ with sets of vertices $\{x,p,q\}$ and
$\{x,q,\ga x\}$.  By comparison, the angles of $\overline{Q}$ at the
vertices $\overline{p}$ and $\overline{q}$ are at least $\frac\pi
2$. If we adjust these angles to $\frac\pi 2$, keeping the lengths of
the three sides $[\overline{x},\overline{p}]$,
$[\overline{p},\overline{q}]$ and $[\overline{q}, \overline{\ga x}]$
fixed, we obtain a quadrilateral $\overline{Q}'$ where the side that
is not adjacent to the right angles has length less than $d(x,\ga x)$.
This gives the lower bound since the length of the side in question is
given by Equation \eqref{eq:planetranslation}.  \cqfd

\medskip The proof of the following result is then similar to that
of Theorem \ref{theo:genHuberabstrait}.

\bcoro \label{coro:encadreloxo}
Let $\wt M$ be a complete simply connected Riemannian manifold
with pinched negative sectional curvature, let $x_0\in\wt M$ and let
$\Ga$ be a nonelementary discrete group of isometries of $\wt M$.
Assume that the Bowen-Margulis measure of $\Ga$ is finite and mixing
for the geodesic flow on $T^1M$.  Let $\KKK$ be a conjugacy class of
loxodromic elements of $\Ga$ with translation length $\ell$.  Then, for
every $\epsilon>0$, if $t$ is big enough,
$$
\frac{i_{\KKK}\,\|\mu_{x_0}\|\,\|\sigma_{\KKK}\|}{\delta_\Ga\,\|m_{\rm BM}
\|\,e^{\frac{\delta_\Ga\,\ell}{2}}} \,e^{\frac{\delta_\Ga}{2}\, t}\,(1-\epsilon) 
\leq
N_{\KKK,\,x_0}(t)\leq\frac{i_{\KKK}\,\|\mu_{x_0}\|\,\|\sigma_{\KKK}\|}
{\delta_\Ga\,\|m_{\rm BM}\|\,(\sinh\frac{\ell}{2})^{\delta_\Ga}}
\,e^{\frac{\delta_\Ga}{2}\, t}\,(1+\epsilon)\,. \;\;\;\Box
$$
\ecoro

In particular, under the assumptions of this result, we have
$$
\lim_{t\ra+\infty} \frac{1}{t}\log N_{\KKK,\,x_0}(t)=\frac{\delta_\Ga}{2}\;.
$$
Theorem \ref{theo:loggrowthintro} in the introduction follows from
this, since if $M=\Ga\bs\wt M$ is a compact manifold, then any
nontrivial element in $\Ga$ is loxodromic, and, as recalled in Section
\ref{sect:counting}, the critical exponent $\delta_\Ga$ is equal to
the topological entropy $h$ of the geodesic flow on $M$, and $m_{\rm BM}$
is finite and mixing.

\medskip
\rem With the notation and definitions of \cite[\S 3.1]{PauPolSha}, if
$\wt F:T^1\wt M\ra \RR$ is a potential (that is, a $\Ga$-invariant
H\"older-continuous map), since the geodesic segment from $x_0$ to
$\alpha x_0$ passes at distance uniformly bounded (by a constant
$c_\ell$ depending only on $\ell$) from the translation axis
$C_\alpha$ of $\alpha$, with $p_\alpha$ the closest point to $x_0$ on
$C_\alpha$, the absolute value of the difference $\int_{x_0}^{\alpha
  x_0}\wt F-\int_{x_0}^{p_\alpha}\wt F-\int_{x_0}^{p_\alpha}\wt
F\circ\iota$ is uniformly bounded (by a constant depending only on
$c_\ell$ and on the maximum of $\wt F$ on the neighbourhood of
$C_\alpha$ of radius $c_\ell$). Hence using the version with potential
of Theorem \ref{theo:mainequicountdown} in \cite[Coro.~20]{ParPau14}
for $\wt F$ and $\wt F\circ\iota$, we have upper and lower bounds for
the asymptotic of the counting function with weigths defined by the
potential: Assume that the critical exponent $\delta_{\Ga,\,F}$ of
$\Ga$ for the potential $\wt F$ is finite and that the Gibbs measure
of $\Ga$ for the potential $\wt F$ is finite and mixing for the
geodesic flow on $T^1M$, then there exists $c>0$ such that for all
$t\geq 0$,
$$
\frac{1}{c}\,e^{\frac{\delta_{\Ga,\,F}}{2}\, t}\leq 
\sum_{\alpha\in\KKK,\; d(x_0,\,\alpha x_0)\leq t} 
m_\alpha\; e^{\int_{x_0}^{\alpha x_0}\wt F} 
\leq c\,e^{\frac{\delta_{\Ga,\,F}}{2}\, t}\;.
$$

\medskip Let us now consider the higher dimensional real hyperbolic
spaces. If $X=\htr$, if $\ga$ is orientation-preserving, and if
$x\in\htr$ is at a distance $s$ from the translation axis of $\ga$,
then
\begin{equation}\label{eq:spacetranslation}
\sinh^2\frac{d(x,\ga x)}2  =
\frac{\sinh^2s\;|e^\lambda-1|^2}{4\,e^\ell}+\sinh^2(\frac\ell 2)\,,
\end{equation}
where $\lambda=\lambda_{\ga}$ is the {\em complex translation length}
of $\ga$, defined as follows.  The loxodromic isometry $\ga$ is
conjugated in the upper halfspace model $\CC\times \mathopen{]}0,
+\infty[$ of $\htr$ to a transformation $(z,r)\mapsto e^{\ell}
(e^{i\theta} z,r)$, where $\theta=\theta_{\ga}\in\RR$ is uniquely
defined modulo $2\pi$, and we define $\lambda= \ell+i\theta\in
\mathopen{]}0, +\infty\mathclose[+i\,\RR/2\pi\ZZ$. Equation
\eqref{eq:spacetranslation} follows from the distance formula in
\cite[pp.~37]{Fenchel89}.

Let $n\in\NN-\{0,1\}$. A loxodromic isometry $\ga$ of $\hnr$ is {\em
  uniformly rotating} if $\ga$ rotates all normal vectors to the
translation axis of $\ga$ by the same angle, called the {\em rotation
  angle} of $\ga$ (which is $0$ if and only if $\ga$ induces the
parallel transport along its translation axis). This property is
invariant under conjugation.

Clearly, all loxodromic isometries of $\hdr$, all
orientation-preserving loxodromic isometries $\htr$, and the
loxodromic isometries of any $\hnr$ with a trivial rotational part,
are uniformly rotating. The orientation-reversing loxodromic
isometries of $\htr$ are not uniformly rotating. More generally, by
the normal form up to conjugation of the elements of $\operatorname{O}
(n-1)$, uniformly rotating orientation-preserving loxodromic
isometries with a nontrivial rotation angle exist in $\hnr$ if and
only if $n$ is odd, and uniformly rotating orientation-reversing
loxodromic isometries exist in $\hnr$ if and only if $n$ is even.  For
a fixed translation length and rotation angle
$\theta\in(\RR-2\pi\ZZ)/(2\pi\ZZ)\,$, with $\theta=\pi$ in the
orientation-reversing case, these elements form a unique conjugacy
class.

Let $\ga$ be a uniformly rotating loxodromic isometry of $\hnr$. Any
configuration that consists of the translation axis of $\ga$, a
geodesic line $L$ orthogonal to the axis and its image $\ga L$ is
contained in an isometrically embedded $\ga$-invariant copy of $\htr$
in $\hnr$ (unique if the rotation angle of $\ga$ is nonzero modulo
$\pi\ZZ$). We then define the {\em complex translation length} of
$\ga$ as the complex translation length of the restriction of $\ga$ to
this subspace.

\blemm\label{lem:offaxis} A uniformly rotating loxodromic isometry
$\ga$ of $\hnr$ with complex translation length $\lambda=
\ell+i\theta$ is $\psi$-equitranslating with $$\psi(t)=\frac{1}{2}
(t-\log(\frac{\cosh\ell-\cos\theta}{2})) +\bigO(e^{-t})$$ as
$t\ra+\infty$.  
\elemm

\dem Let $x$ be a point in $\hnr$ at a distance $s$ from the
translation axis of $\ga$. We only have to prove that, as
$s\to+\infty$,
$$
d(x,\ga x)= 2s+\log\big(\frac{\cosh\ell-\cos\theta}{2}\big)
+\bigO(e^{-2s})\;.
$$
As noted above, it suffices to consider the case $n=3$. By Equation
\eqref{eq:spacetranslation}, we have
$$
\frac{e^{d(x,\,\ga x)}}4 =
\frac{e^{2\,s}\,(e^{2\ell}-2e^\ell\cos\theta+1)}{16\,e^\ell}+
\bigO(1)\,,
$$
as $s\to+\infty$, which proves the claim after simplification and
taking the logarithm.  
\cqfd

\bcoro \label{coro:mainloxo} Let $\Ga$ be a nonelementary discrete
group of isometries of $\hnr$, whose Bowen-Margulis measure is finite,
and let $x_0\in\hnr$.  Let $\KKK$ be a conjugacy class of uniformly
rotating loxodromic elements of $\Ga$ with complex translation length
$\lambda= \ell+i\theta$.  Then, as $t\ra+\infty$,
$$
N_{\KKK,\,x_0}(t)\sim 
\frac{2^{\frac{\delta_\Ga}2}\,i_\KKK\,\|\mu_{x_0}\|\,\|\sigma_{\KKK}\|}
{\delta_\Ga\,\|m_{\rm BM}\|\,
(\cosh\ell-\cos\theta)^{\frac{\delta_\Ga}{2}}}
\;e^{\frac{\delta_\Ga}{2}\, t}\,.
$$
If $\Ga$ is arithmetic or if $M$ is compact, then the error term is
$\bigO(e^{(\frac{\delta_\Ga}{2} -\kappa) t})$ for some $\kappa >0$.
Furthermore, if $v_\alpha$ is the unit tangent vector at $x_0$ to the
geodesic segment $[x_0,\alpha x_0]$ for every $\alpha\in \Ga-
\Ga_{x_0}$, for the weak-star convergence of measures on $T^1\wt M$,
we have
$$
\lim_{t\ra+\infty}\; \frac{\delta_\Ga\;\|m_{\rm BM}\|\,
(\cosh\ell-\cos\theta)^{\frac{\delta_\Ga}{2}}}
{2^{\frac{\delta_\Ga}2}\,i_{\KKK}\,\|\sigma_{\KKK}\|}
\;e^{-\frac{\delta_\Ga}{2}\, t}\;
\sum_{\alpha\in\KKK,\; 0<d(x_0,\,\alpha x_0)\leq t} 
m_\alpha \,\Delta_{v_\alpha}\;=\; (\pi_+^{-1})_*\mu_{x_0}\;.
$$
\ecoro 

\dem As mentioned in Section \ref{sec:rappels}, since $\hnr$ has
constant sectional curvature, the Bowen-Margulis measure of $\Ga$,
since finite, is mixing for the geodesic flow on $T^1M$. We have
already seen that $i_\KKK$ is finite and that $\|\sigma_\KKK\|$ is
positive and finite. The result follows from Theorem
\ref{theo:genHuberabstrait} and Lemma \ref{lem:offaxis}.  \cqfd

\medskip \rem Let $\Ga$ be a group of isometries of $X$ and assume
that $\ga$ is a loxodromic element of $\Ga$. The element $\ga$ is
$\Ga${\em-reciprocal} if there exists an element in $\Ga$ that
switches the two fixed points of $\ga$.  If $\ga$ is reciprocal, then
let $\iota_\Ga(\ga)=2$, otherwise, we set $\iota_\Ga(\ga)=1$. The
stabiliser in $\Ga$ of the translation axis $C_{\ga}$ of $\ga$ is
generated by the maximal cyclic subgroup of $\Ga$ containing $\ga$, by
an elliptic element that switches the two points at infinity of
$C_{\ga}$ if $\ga$ is $\Ga$-reciprocal, and a (possibly trivial) group
of finite order, which is the pointwise stabiliser
$\operatorname{Fix}_\Ga (C_\ga)$ of $C_{\ga}$. Thus, if $\KKK$ is the
conjugacy class of $\ga$ in $\Ga$,
$$
\iota_\KKK=\iota_\Ga(\ga)
[\operatorname{Fix}_\Ga (C_\ga):
\operatorname{Fix}_\Ga (C_\ga)\cap Z_{\Ga}(\ga)]\,.
$$
In particular, if $n=2$, or if $n=3$ and $\ga$ preserves the
orientation, then $\iota_\KKK=\iota_\Ga(\ga)$. Hence Theorem
\ref{theo:intro} in the Introduction when $\ga_0$ is loxodromic
follows from Corollary \ref{coro:mainloxo}.

\medskip When $\Ga$ has finite covolume, the constant in Corollary
\ref{coro:mainloxo} can be made more explicit.

\bcoro\label{corogenHuber} Let $\Ga$ be a discrete group of isometries
of $\hnr$ with finite covolume and let $x_0\in\hnr$.  Let $\KKK$ be
the conjugacy class of a uniformly rotating loxodromic element $\ga_0$
of $\Ga$ with complex translation length $\lambda=\ell+i\theta$, let
$m_{\ga_0}$ be the order of $\ga_0$ in the maximal cyclic group
containing $\ga_0$, and let $n_{\ga_0}$ be the order of the
intersection of the pointwise stabiliser of the translation axis of
$\ga_0$ with the centraliser of $\ga_0$. Then, as $t\ra+\infty$,
$$
N_{\KKK,\,x_0}(t)\sim \frac{\Vol(\SSS^{n-2})\;\ell}
{2^{\frac{n-1}2}\,(n-1)\,m_{\ga_0}\,n_{\ga_0}\,\Vol(\Ga\bs\hnr)\,
(\cosh \ell-\cos \theta)^{\frac{n-1}2}}\;e^{\frac{n-1}2\, t }\,.
$$
If $\Ga$ is arithmetic or if $M$ is compact, then the error term is
$\bigO(e^{(\frac{n-1}{2} -\kappa) t})$ for some $\kappa >0$.
Furthermore, if $\Ga_{x_0}$ is trivial, if $v_\alpha$ is the unit
tangent vector at $x_0$ to the geodesic segment $[x_0,\alpha x_0]$ for
every $\alpha\in \Ga- \{e\}$, with $\Vol_{T^1_{x_0}\hnr}$ the
spherical measure on $T^1_{x_0}\hnr$, we have, for the weak-star
convergence of measures on $T^1_{x_0}\hnr$,
\begin{multline*}
\frac{(n-1)\,m_{\ga_0}\,n_{\ga_0}\,\Vol(\SSS^{n-1})\,\Vol(\Ga\bs\hnr)\,
(\cosh\ell-\cos\theta)^{\frac{n-1}{2}}}
{2^{\frac{1-n}2}\,\Vol(\SSS^{n-2})\,\ell\;e^{\frac{n-1}{2}\, t}}
\sum_{\alpha\in\KKK,\; 0<d(x_0,\,\alpha x_0)\leq t} \Delta_{v_\alpha}
\\\stackrel{*}{\rightharpoonup}\;\; \Vol_{T^1_{x_0}\hnr}\;.
\end{multline*}
\ecoro

\dem Since $\Ga$ has finite covolume, we have $\delta_\Ga=n-1$ and we
can normalise the Patterson-Sullivan measure $\mu_{x_0}$ at $x_0$ to
have total mass $\Vol(\SSS^{n-1})$, so that $(\pi_+^{-1})_*\mu_{x_0}=
\Vol_{T^1_{x_0}\hnr}$.  By \cite[Prop.~10, 11]{ParPauRev}, we have
$$
\|m_{\rm BM}\|=2^{n-1}\Vol(\SSS^{n-1})\Vol(\Ga\bs\hnr)
$$ 
and
$$
\|\sigma_{C_{\ga_0}}\|= \Vol(\SSS^{n-2}) \frac
\ell{|\operatorname{Fix}_{\Ga_0} (C_{\ga_0})|\,\iota_\Ga(\ga_0)\,m_{\ga_0}}\,,
$$ 
since $\Vol(\Ga_{C_{\ga_0}}\bs C_{\ga_0}) =\frac{\ell}
{\iota_\Ga(\ga_0) \,m_{\ga_0}}$.  The claims hence follow from the
previous remark and from Corollary \ref{coro:mainloxo}.  \cqfd

\medskip The proof of the loxodromic case of Corollary
\ref{coro:intro} of the Introduction follows from Corollary
\ref{corogenHuber} by taking $n=2$, $\Ga$ torsion-free (so that
$n_{\ga_0}=1$), and $\ga_0$ primitive (so that $m_{\ga_0}=1$) and
orientation-preserving (so that $\cos\theta=1$). The area of a
complete, connected, finite area hyperbolic surface with genus $g$ and
$p$ punctures is $2\pi(2g+p-2)$.

%In the case of compact surfaces, Huber \cite{Huber98}  gives an 
%explicit error bound of the form $\bigO(e^{\frac{3}{8} t})$ .

\section{The geometry of parabolic isometries}
\label{sect:parabolic}

In this section, we fix a parabolic isometry $\ga$ of a complete
$\CAT(-1)$ geodesic metric space $X$. We fix a horoball $C_\ga$
centred at the fixed point of $\ga$, and we call {\em horospherical
  translation length} of $\ga$ the quantity
$$
\ell=\ell_{\ga}=\inf_{y\in \partial C_\ga} d(y, \ga y)\;.
$$

We will say that $\ga$ is {\em uniformly translating} if $d(y,\ga y)$
is independent of $y\in \partial C_\ga$. Note that being uniformly
translating does not depend on the choice of $C_\ga$, 
but the value of
$\ell$ does (and can be fixed arbitrarily in $]0,+\infty[$ when $X$ is
a Riemannian manifold). 

Every parabolic isometry of $X=\hdr,\htr$ is uniformly translating, but using Euclidean screw motions,
there exist parabolic isometries in $X=\HH^4_\RR$ which are not
uniformly translating (and the map $y\mapsto d(y, \ga y)$ is not even
bounded). If $X=\hnr$ and  if $\ga$ induces a Euclidean translation on
$\partial C_\ga$, then $\ga$ is uniformly translating. Recall that by
Bieberbach's theorem, any discrete group of isometries of $\hnr$,
preserving a given horosphere and acting cocompactly on it, contains a
finite index subgroup consisting of uniformly
translating parabolic isometries and the identity.

If $X=\hdr$, if $x\in X$ is at a distance $s$ from the horoball
$C_\ga$, then
\begin{equation}\label{eq:planeparabolic}
d(x,\ga x)=2\arsinh(e^s\,\sinh\frac{\ell} 2)\,.
\end{equation}
This is immediate by considering the upper halfplane model and
assuming that $\ga$ has fixed point $\infty$,  by applying twice
\cite[Thm.~7.2.1 (iii)]{Beardon83}. A similar triangle inequality and
comparison argument as in the proof of Lemma
\ref{lem:generaltranslation} shows the following result.

\blemm\label{lem:generalparabolic} If $x\in X$ is at distance $s>0$
from the horoball $C_\ga$, if $p_\ga$ is the closest point to $x$ on
$C_\ga$, then
$$
2\arsinh(e^s\sinh\frac{\ell}  2)\le d(x,\ga x)
\le 2s+ d(p_\ga,\ga \,p_\ga)\,. \;\;\Box
$$
\elemm

\bcoro\label{lem:offhorosphere} A uniformly translating parabolic
isometry $\ga$ of $\hnr$ with horospherical translation length $\ell$
is $\psi$-equitranslating with $$\psi(t)= \frac{t}{2}-
\log(\sinh\frac{\ell}{2}) -\log 2 +\bigO(e^{-t})$$ as
$t\ra+\infty$.  
\ecoro

\dem Let $x$ be a point in $\hnr$ at a distance $s$ from
the horoball $C_\ga$. We only have to prove that, as
$s\to+\infty$,
$$
d(x,\ga x)= 2s+2\log\big(\sinh\frac{\ell}{2}\big)+2\log 2
+\bigO(e^{-2s})\;.
$$
It suffices to consider the case $n=2$ (the points $x,\ga x$ and the
fixed point of $\ga$ are contained in a copy of $\hdr$), in which case
the result follows from Equation \eqref{eq:planeparabolic}. \cqfd

\medskip
\rem 
If $X=\wt M$ and $\Ga$ are as in Section \ref{sect:counting}, if $\ga$
is a parabolic isometry of $\Ga$ and if $\KKK$ is the conjugacy class
of $\ga$ in $\Ga$, the quantities $\|\sigma_{\KKK}\|$ and $i_\KKK$
defined in that Section are not always finite.  Note that
$\|\sigma_{\KKK}\|$ is positive, since $\Ga$ is nonelementary.

$\bullet$~ The mass $\|\sigma_{\KKK}\|$ is finite for instance if the
fixed point $\xi_\ga$ of $\ga$ is a bounded parabolic fixed point
(that is, if its stabiliser $\Ga_{\xi_\ga}$ in $\Ga$ acts cocompactly
on $\Lambda\Ga-\{\xi_\ga\}$), which is in particular the case if $\Ga$
is a lattice or is geometrically finite.

$\bullet$~ The index $i_\KKK$ is equal to $1$ if $\ga$ is central in
the stabiliser $\Ga_{C_\ga}$ of the horoball $C_\ga$. This is in
particular the case, up to passing to a finite index subgroup of
$\Ga$, if $\Ga$ is a lattice or is geometrically finite, as well as if
$X$ is a symmetric space and $\ga$ is in the center of the nilpotent
Lie group of isometries of $X$ acting simply transitively on the
horosphere $C_\ga$ (see Proposition \ref{prop:heisenberg} below: in
the complex hyperbolic space $\hnc$, this center consists of the
vertical Heisenberg translations). If $X=\hdr$, we have $i_\KKK=1$ if
no nontrivial elliptic element of $\Ga$ fixes $\xi_\ga$ (in particular
if $\Ga$ is torsion-free), and $i_\KKK=2$ otherwise. In the complex
hyperbolic space $\hnc$, the stabilisers of horoballs are not abelian
and $i_\KKK$ is finite only if $\KKK$ consists of vertical Heisenberg
translations.

\medskip
A proof similar to that of Corollary \ref{coro:mainloxo}
gives the following result, which implies in particular Theorem
\ref{theo:intro} in the Introduction when $\ga_0$ is parabolic.

\bcoro \label{coro:mainpara} Let $\Ga$ be a nonelementary discrete
group of isometries of $\hnr$, whose Bowen-Margulis measure is finite,
and let $x_0\in\hnr$.  Let $\KKK$ be a conjugacy class of uniformly
translating parabolic elements of $\Ga$ with horospherical
translation length $\ell$, with $\|\sigma_{\KKK}\|$ and $i_\KKK$
finite.  Then, as $t\ra+\infty$,
$$
N_{\KKK,\,x_0}(t)\sim 
\frac{i_\KKK\,\|\mu_{x_0}\|\,\|\sigma_{\KKK}\|}
{\delta_\Ga\,\|m_{\rm BM}\|\,
(2\sinh\frac{\ell}{2})^{\delta_\Ga}}
\;e^{\frac{\delta_\Ga}{2}\, t}\,.
$$
If $\Ga$ is arithmetic, then the error term is
$\bigO(e^{(\frac{\delta_\Ga}{2} -\kappa) t})$ for some $\kappa
>0$. Furthermore, if $v_\alpha$ is the unit tangent vector at $x_0$ to
the geodesic segment $[x_0,\alpha x_0]$ for every $\alpha\in \Ga-
\Ga_{x_0}$, for the weak-star convergence of measures on $T^1\wt M$,
we have
$$
\lim_{t\ra+\infty}\; \frac{\delta_\Ga\;\|m_{\rm BM}\|\,
(2\,\sinh\frac{\ell}{2})^{\delta_\Ga}}
{i_{\KKK}\,\|\sigma_{\KKK}\|} \;e^{-\frac{\delta_\Ga}{2}\, t}\;
\sum_{\alpha\in\KKK,\; 0<d(x_0,\,\alpha x_0)\leq t} 
m_\alpha\,\Delta_{v_\alpha}
\;=\; (\pi_+^{-1})_*\mu_{x_0}\;. \;\;\Box
$$
\ecoro

\bcoro\label{coro:finitevolpara} Let $\Ga$ be a discrete group of
isometries of $\hnr$ with finite covolume and let $x_0\in\hnr$.  Let
$\KKK$ be the conjugacy class of a uniformly translating parabolic
element $\ga_0$ of $\Ga$ with $i_\KKK$ finite. Then, as $t\ra+\infty$,
$$
N_{\KKK,\,x_0}(t) \sim
\frac{i_\KKK\,\Vol(\Ga_{C_{\ga_0}}\bs C_{\ga_0})}
{\Vol(\Ga\bs\hnr)\,(2\sinh\frac{\ell}{2})^{n-1}}\;e^{\frac{n-1}2\, t }\,.
$$
If $\Ga$ is arithmetic, then the error term is $\bigO(
e^{(\frac{n-1}{2} -\kappa) t})$ for some $\kappa >0$.  Furthermore, if
$\Ga_{x_0}$ is trivial, if $v_\alpha$ is the unit tangent vector at
$x_0$ to the geodesic segment $[x_0,\alpha x_0]$ for every $\alpha\in
\Ga- \{e\}$, with $\Vol_{T^1_{x_0}\hnr}$ the spherical measure on
$T^1_{x_0}\hnr$, we have, for the weak-star convergence of measures on
$T^1_{x_0}\hnr$,
$$
\frac{\Vol(\SSS^{n-1})\,\Vol(\Ga\bs\hnr)\,
(2\sinh\frac{\ell}{2})^{\frac{n-1}{2}}}
{i_\KKK\,\Vol(\Ga_{C_{\ga_0}}\bs C_{\ga_0})}
\;e^{-\frac{n-1}{2}\, t}
\sum_{\alpha\in\KKK,\; 0<d(x_0,\,\alpha x_0)\leq t} \Delta_{v_\alpha}
\stackrel{*}{\rightharpoonup}\;\; \Vol_{T^1_{x_0}\hnr}\;.
$$
\ecoro

\dem The claims  are proved in the same way as Corollary
\ref{corogenHuber}, using the equality $$\|\sigma_{\KKK}\|=
2^{n-1}\,(n-1)\,\Vol(\Ga_{C_{\ga_0}}\bs C_{\ga_0})\,,$$ see
\cite[Prop.~29 (2)]{ParPau14}.  
\cqfd

\medskip 
The parabolic case of Corollary \ref{coro:intro} of the Introduction
follows from Corollary \ref{coro:finitevolpara}.  Consider the upper
halfplane model of $\hdr$ and normalise the group such that $\ga_0$ is
the translation $z\mapsto z+1$. We choose $C_{\ga_0}$ to be the
horoball that consists of points with imaginary part at least
$1$. Since $\ga_0$ is primitive and $\Ga$ is torsion-free, we have
$\Ga_{C_{\ga_0}}={\ga_0}^\ZZ$ and $i_\KKK=1$. Hence $\Vol(\Ga_{C_{\ga_0}}\bs
C_{\ga_0})=1$ by a standard computation of hyperbolic area. Now,
$\sinh \frac{\ell}{2}= \frac{1}{2}$, and the claim follows as in the
proof of the loxodromic case after Corollary \ref{corogenHuber}.

\medskip 
We end this section by giving a necessary and sufficient criterion for
a parabolic isometry of the complex hyperbolic space $\hnc$ to be
uniformly translating. We refer to \cite{Goldman99}, besides the
reminder below, for the basic properties of $\hnc$.

On $\CC^{n+1}=\CC\times\CC^{n-1}\times\CC$, consider the Hermitian
product with signature $(1,n)$ defined by
$$
\langle z,w\rangle=-z_0\ov w_n+z\cdot\ov w-z_n\ov w_0\,,
$$
where $(z,w)\mapsto z\cdot\ov w$ is the standard Hermitian scalar
product on $\CC^{n-1}$. Let $q(z)=\langle z, z\rangle$ be the
corresponding Hermitian form.  The projective model of the complex
hyperbolic space $\hnc$ corresponding to this choice of $q$ is the set
$$
\{[w_0:w:1]\in\PP_n(\CC)\;:\; q(w_0,w,1)<0\}\,,
$$
endowed with the Riemannian metric, normalised to have sectional
curvature between $-4$ and $-1$, whose Riemannian distance is given by
$$
d(X,Y)=\arcosh\sqrt{\frac{\langle x,y\rangle\langle y,x\rangle}
{q(x)\,q(y)}}
$$
for any representatives $x,y$ of $X,Y$ in $\CC^{n+1}$, see
\cite[p.~77]{Goldman99}, where the sectional curvature is normalised
to be between $-1$ and $-\frac 14$.  The boundary at infinity of
$\hnc$ is
$$
\partial_\infty\hnc=
\{[w_0:w:1]\in\PP_n(\CC)\;:\; q(w_0,w,1)=0\}\cup\{\infty\}\,,
$$
where $\infty=[1:0:0]$. For every $s>0$, the set 
$$
\H_s=\{[w_0:w:1]\in\PP_n(\CC)\;:\; q(w_0:w:1)=-s\}
$$
is a horosphere centred at $\infty$. 

The parabolic isometries $\ga$ of $\hnc$ fixing $\infty$ are the
mappings induced by the projective action of the matrices
\begin{equation}\label{eq:para}
\wt \ga=\begin{pmatrix} 1& \;a^*&z_0\\0&A&b\\0&0&1\end{pmatrix}\,,
\end{equation}
where $A\in U(n-1)$, $a^*=\;^t\overline{a}$ and $A a= b$, see \cite[\S
4.1]{CheGre74} and \cite[p.~371]{ParPau10GT}. For every $Z=[z_0:z:1]\in \partial_\infty \hnc -
\{\infty\}$, the isometry induced by the matrix
$$
T_Z=\begin{pmatrix} 1&\;z^*&z_0\\0&1&z\\0&0&1
\end{pmatrix}
$$
is called a {\it Heisenberg translation}, which is {\it
  vertical} if $z=0$.  The group of Heisenberg translations (which
identifies with the Heisenberg group of dimension $2n-1$, see
\cite{Goldman99}) acts simply transitively on $\partial_\infty\hnc
-\{\infty\}$ and on each horosphere $\H_s$ for $s>0$.

\bprop\label{prop:heisenberg} 
A parabolic isometry $\ga$ of the complex hyperbolic space $\hnc$ is
uniformly translating if and only if it is a vertical Heisenberg
translation. Furthermore, if $\ga$ is not a vertical Heisenberg
translation, then the map $y\mapsto d(y,\ga y)$ is unbounded on any
horosphere of $\hnc$ centred at the fixed point of $\ga$.  
\eprop

\dem 
For all $W=[w_0:w:1]\in\H_2$ and any parabolic isometry $\ga$ as given
by Equation \eqref{eq:para}, we have
$$
d(W,\ga W)=
%\arcosh\frac{|w^*A^*w-|w|^2+\bigO(|w|)|}2\,.
\arcosh\frac{|w^*(A^*-I)w+\bigO(|w|)|}2\,.
$$
If $A$ is not the identity, then $w^*(A^*-I)w$ is equivalent to $|w|^2$ (up
to a positive constant) on some line in $\CC^{n-1}$, which makes the map
$W\mapsto d(W,\ga W)$ unbounded on $\H_2$.  Thus we are reduced to
considering Heisenberg translations. For all $Z=[z_0:z:1]
\in \partial_\infty \hnc - \{\infty\}$, we have
$$
d(W,T_Z W)=\arcosh\frac{|z\cdot\ov w-\ov zw-z_0-2|}2\,.
$$
It is easy to see that this distance is independent of $W$ if and only
if $z=0$, and is unbounded otherwise.  \cqfd

\section{The geometry of elliptic isometries}
\label{sect:elliptic}

In this section, we fix $n\geq 2$ and a nontrivial elliptic isometry
$\ga$ of $\hnr$. We denote by $C_\ga$ the fixed point set of $\ga$,
which is a nonempty proper totally geodesic subspace of $\hnr$ of
dimension $k=k_\ga$.

We will say that $\ga$ is {\em uniformly rotating} if there exists
$\theta=\theta_\ga\in\mathopen{]}0,\pi]$ (called the {\it rotation
  angle} of $\ga$) such that for every $v\in \normalout C_\ga$, the
(nonoriented) angle between $v$ and $\ga v$ is $\theta$. This property
is invariant under conjugation, and once $k$ and $\theta$ are fixed,
there exists only one conjugacy class of uniformly rotating nontrivial
elliptic isometries. Note that when $n=2$ or $n=3$, every elliptic
isometry $\ga$ is uniformly rotating, and $\theta=\pi$ if $\ga$ does
not preserve the orientation. But there exist elliptic isometries in
$\HH^4_\RR$ which are not uniformly rotating.

Assume that $\ga$ belongs to a nonelementary discrete group of
isometries $\Ga$ of $\hnr$, and let $\KKK$ be the conjugacy class of
$\ga$ in $\Ga$.

$\bullet$~ The skinning measure $\|\sigma_\KKK\|$ is positive if and
only if $\Lambda\Ga$ is not contained in $\partial_\infty C_\ga$. This
is in particular the case if $n=2$.  Furthermore, $\|\sigma_\KKK\|$ is
finite for instance if $\Ga_{C_\ga}\bs C_\ga$ is compact or if
$\partial_\infty C_\ga\cap \Lambda\Ga$ is empty. This is in particular
the case if $n=2$ and if $\ga$ preserves the orientation. But when
$n=2$ and $\ga$ does not preserve the orientation, the measure
$\|\sigma_\KKK\|$ is not necessary finite. 

For instance, let $\Ga= T(\infty,\infty,\infty)$ be the discrete group
of isometries of $\HH^2_\RR$ generated by the reflexions $s_1, s_2,
s_3$ on the sides of an ideal hyperbolic triangle. Then $C_{s_1}$ is
one of these sides. Let us prove that $\pi_*\wt{\sigma}_{C_{s_1}}$ is
a constant multiple of the Lebesgue measure along $C_{s_1}$. Indeed,
the Patterson-Sullivan measure at infinity of the disc model of
$\HH^2_\RR$ based at its origin is a multiple of the Lebesgue measure
$d\theta$ on the circle, since $\Ga$ has finite covolume. Since
$d\theta$ is conformally invariant under every isometry of
$\HH^2_\RR$, the measure $\pi_*\wt{\sigma}_{C_{s_1}}$ on $C_{s_1}$ is
invariant under every loxodromic isometry preserving $C_{s_1}$, hence
the result.  Since $C_{s_1}$ injects in $\Ga\bs\HH^2_\RR$ and since
its stabiliser in $\Ga$ has order $2$, the measure
$\pi_*\sigma_{C_{s_1}}$ is the multiple by half the above constant of
the Lebesgue measure on the image of $C_{s_1}$ in $\Ga\bs\HH^2_\RR$,
which is infinite.

$\bullet$~ If $n=2$ and $k_\ga=1$ (so that $\ga$ reverses the
orientation), then every isometry preserving $C_{\ga}$ commutes with
$\ga$, hence $i_\KKK=1$. If $n=2$ and $k_\ga=0$ (so that $\ga$
preserves the orientation), then the finite group $\Ga_{C_{\ga}}$ is
either cyclic, in which case $\Ga_{C_{\ga}}=Z_\Ga(\ga)$ and
$i_\KKK=1$, or it is dihedral. Assume the second case
holds. If the rotation angle of $\ga$ is $\pi$, then again
$\Ga_{C_{\ga}}=Z_\Ga(\ga)$ and $i_\KKK=1$. Otherwise, $i_\KKK=2$.

\blemm\label{lem:offfixedpoint} 
A uniformly rotating elliptic isometry $\ga$ of $\hnr$ with rotation
angle $\theta$ is $\psi$-equitranslating with
$$
\psi(t)=\frac{t}{2} -\log\frac{\sin\theta}{2} +\bigO(e^{-\frac{t}{2}})
$$ 
as $t\ra+\infty$.  
\elemm

\dem By the formulas in right-angled hyperbolic triangles (see
\cite[Theo.~7.11.2 (ii)]{Beardon83}), if $x\in \hnr$ is at distance
$s$ from the fixed point set $C_\ga$ of $\ga$, we have
$$
\sinh\frac{d(x,\ga x)}{2}=\sinh s\;\frac{\sin\theta}{2}\;.
$$
The result follows as in Lemma \ref{lem:offaxis}.
\cqfd

\medskip The next result follows from this lemma in the same way as
Corollary \ref{coro:mainloxo} follows from Lemma \ref{lem:offaxis}. It
implies Theorem \ref{theo:intro} in the Introduction when $\ga_0$ is
elliptic.

\bcoro \label{coro:mainellip} Let $\Ga$ be a nonelementary discrete
group of isometries of $\hnr$, whose Bowen-Margulis measure is finite,
and let $x_0\in\hnr$.  Let $\KKK$ be a conjugacy class of uniformly
rotating nontrivial elliptic elements of $\Ga$ with rotation angle
$\theta$, such that $\|\sigma_\KKK\|$ and $i_\KKK$ are positive and 
finite.  Then, as $t\ra+\infty$,
$$
N_{\KKK,\,x_0}(t)\sim 
\frac{i_\KKK\,\|\mu_{x_0}\|\,\|\sigma_{\KKK}\|}
{\delta_\Ga\,\|m_{\rm BM}\|\,
(\sin\frac{\theta}{2})^{\delta_\Ga}}
\;e^{\frac{\delta_\Ga}{2}\, t}\,.
$$
If $\Ga$ is arithmetic or if $M$ is compact, then the error term is
$\bigO(e^{(\frac{\delta_\Ga}{2} -\kappa) t})$ for some $\kappa
>0$. Furthermore, if $v_\alpha$ is the unit tangent vector at $x_0$ to
the geodesic segment $[x_0,\alpha x_0]$ for every $\alpha\in \Ga-
\Ga_{x_0}$, for the weak-star convergence of measures on $T^1\wt M$,
we have
$$
\lim_{t\ra+\infty}\; \frac{\delta_\Ga\;\|m_{\rm BM}\|\,
(\sin\frac{\theta}{2})^{\delta_\Ga}}
{i_{\KKK}\,\|\sigma_{\KKK}\|}
\;e^{-\frac{\delta_\Ga}{2}\, t}\;
\sum_{\alpha\in\KKK,\; 0<d(x_0,\,\alpha x_0)\leq t} 
m_\alpha\,\Delta_{v_\alpha}
\;=\; (\pi_+^{-1})_*\mu_{x_0}\;.\;\;\;\Box
$$
\ecoro

\section{Counting conjugacy classes of subgroups}
\label{sect:countingsubgroups}

Let $\wt M, x_0,\Ga$ be as in the beginning of Section
\ref{sec:rappels}. Let $\Ga_0$ be a subgroup of $\Ga$, and let 
$$\KKK=\{\ga\Ga_0\ga^{-1}\;:\;\ga\in\Ga\}$$ be its conjugacy class in
$\Ga$. In this Section, we will study the asymptotic growth, as
$t\ra+\infty$, of the cardinality of 
$$
\{A\in\KKK\;:\; \inf_{\alpha\in
  A -\{e\}} d(x_0,\,\alpha x_0)\leq t\}\;,
$$ 
the set (assumed to be finite) of the conjugates of $\Ga_0$ in $\Ga$
whose minimal displacement of $x_0$ is at most $t$.

We will assume the following conditions on $\Ga_0$: 

($*$)~ There exists a nonempty proper closed convex subset $C_0$ in
$\wt M$ such that the normaliser $N_\Ga(\Ga_0)$ of $\Ga_0$ in $\Ga$ is
a subgroup of the stabiliser $\Ga_{C_0}$ of $C_0$ in $\Ga$, with
finite index, denoted by $i_0$, and such that the family $(\ga
C_0)_{\ga\in\Ga/\Ga_{C_0}}$ is locally finite in $\wt M$;

\smallskip 
($**$)~ There are $c_-,c_+\in\mathopen{]}0,+\infty[$ such
that $c_-\leq \inf_{\ga\in \Ga_0 -\{e\}} d(y,\,\ga y)\leq c_+$ for
every $y\in \partial C_0$.

\medskip 
For instance, $\Ga_0$ could be an infinite index malnormal
torsion-free cocompact stabiliser of a proper totally geodesic
subspace $C_0$ of dimension at least $1$ in $\wt M$, or a torsion-free
cocompact stabiliser of a horosphere centered at a parabolic fixed
point of $\Ga$ (with $C_0$ the horoball bounded by this horosphere),
in which cases $i_0=1$ and $ \|\sigma_{C_0}\|$ is positive and finite.

For every $A=\ga\Ga_0\ga^{-1}\in\KKK$, let 
$$
m_A= (\card(\Ga_{x_0}\cap\Ga_{\ga C_0}))^{-1}\,,
$$ 
which is well-defined since the normaliser of $\Ga_0$ in $\Ga$
stabilises $C_0$. We define the counting function
$$
N_{\KKK,\,x_0}(t)=
\sum_{A\in\KKK,\; \inf_{\alpha\in A -\{e\}} d(x_0,\,\alpha x_0)\leq t} m_A
=\sum_{\ga\in\Ga/N_\Ga(\Ga_0),\; \inf_{\alpha\in \Ga_0 -\{e\}} 
d(x_0,\,\ga\alpha\ga^{-1} x_0)\leq t} m_A\,.
$$

\bprop Let $\wt M$ be a complete simply connected Riemannian manifold
with pinched negative sectional curvature, let $x_0\in\wt M$,  and let
$\Ga$ be a nonelementary discrete group of isometries of $\wt M$.
Assume that the Bowen-Margulis measure of $\Ga$ is finite and mixing
for the geodesic flow on $T^1M$.  Let $\Ga_0$ be a subgroup of $\Ga$
and let $C_0$ be a subset of $\wt M$ satisfying the conditions ($*$) and
($**$), such that the skinning measure $\|\sigma_{C_0}\|$ is positive
and finite. Let $\KKK$ be the conjugacy class of $\Ga_0$ in $\Ga$.
Then, for every $\epsilon>0$, if $t$ is big enough,
$$
\frac{i_{0}\,\|\mu_{x_0}\|\,\|\sigma_{C_0}\|}{\delta_\Ga\,\|m_{\rm BM}
\|\,e^{\frac{\delta_\Ga\,c_+}{2}}} \,e^{\frac{\delta_\Ga}{2}\, t}\,(1-\epsilon) 
\leq
N_{\KKK,\,x_0}(t)\leq\frac{i_{0}\,\|\mu_{x_0}\|\,\|\sigma_{C_0}\|}
{\delta_\Ga\,\|m_{\rm BM}\|\,(\sinh\frac{c_-}{2})^{\delta_\Ga}}
\,e^{\frac{\delta_\Ga}{2}\, t}\,(1+\epsilon)\,. 
$$
\eprop

\dem Let $\ga\in\Ga$. By the local finiteness assumption, except for
finitely many cosets of $\ga$ in $\Ga/\Ga_{C_0}$, the point $x_0$ does
not belong to $\ga C_0$. As in Lemma \ref{lem:generaltranslation}, if
$x_0\in X$ is at distance $s$ from $\ga C_0$, we have
$$
2\arsinh(\cosh s\sinh\frac{c_-}2)\leq
\inf_{\alpha\in \Ga_0 -\{e\}} d(x_0,\ga\alpha\ga^{-1} x_0)
\leq 2s+c_+\,.
$$
The proof is then similar to the proof of Corollary
\ref{coro:encadreloxo}. \cqfd

\medskip We have the following more precise result under stronger
assumptions on $\Ga_0$, with a proof similar to those of Corollaries
\ref{coro:mainloxo} and \ref{coro:mainpara}.

\btheo Let $\Ga$ be a nonelementary discrete group of isometries of
$\hnr$ with finite Bowen-Margulis measure, and let $x_0\in\hnr$.  Let
$\Ga_0$ be the stabiliser in $\Ga$ of a bounded parabolic fixed point
of $\Ga$, acting purely by translations on the boundary of any
horoball $C_0$ centred at this fixed point. Let $\KKK$ be the
conjugacy class of $\Ga_0$ in $\Ga$ and let $\ell=\min_{\ga\in\Ga_0
  -\{e\}} d(y,\ga y)$ for any $y\in \partial C_0$.  Then, as
$t\ra+\infty$,
$$
N_{\KKK,\,x_0}(t)\sim \frac{\|\mu_{x_0}\|\,\|\sigma_{C_0}\|}
{\delta_\Ga\,\|m_{\rm BM}\|\,
(2\sinh\frac{\ell}{2})^{\delta_\Ga}}
\,e^{\frac{\delta_\Ga}{2}\, t}\,.
$$
If $\Ga$ is arithmetic, then the error term is
$\bigO(e^{(\frac{\delta_\Ga}{2} -\kappa) t})$ for some $\kappa
>0$. \cqfd 
\etheo

%{\small\bibliography{../biblio} }
{\small \bibliography{../viitteet} }

\bigskip
{\small
\noindent \begin{tabular}{l} 
Department of Mathematics and Statistics, P.O. Box 35\\ 
40014 University of Jyv\"askyl\"a, FINLAND.\\
{\it e-mail: jouni.t.parkkonen@jyu.fi}
\end{tabular}
\medskip

\noindent \begin{tabular}{l}
D\'epartement de math\'ematique, UMR 8628 CNRS, B\^at.~425\\
Universit\'e Paris-Sud,
91405 ORSAY Cedex, FRANCE\\
{\it e-mail: frederic.paulin@math.u-psud.fr}
\end{tabular}
}

\end{document}